\documentclass[11pt]
{article}
\usepackage{amsthm,amsfonts,amssymb, latexsym}

\newcommand{\n}{\noindent}
\newcommand{\vn}{\vspace{5mm} \noindent}

\def\Aut{{\rm Aut}}

\def\Hom{{\rm Hom}}
\def\End{{\rm End}}

\def\Spec{{\rm Spec}}
\def\Spf{{\rm Spf}}
\def\Isom{{\rm Isom}}

\def\dim{{\rm dim}}
\def\rk{{\rm rk}}
\def\sdim{{\rm sdim}}
\def\Def{{\rm Def}}

\def\Ima{{\rm Im}}
\def\csd{{\rm completely slope divisible}}
\def\gfc{{\rm geometrically fiberwise constant}}
\def\rk{{\rm rk}}
\def\deg{{\rm deg}}
\def\Ker{{\rm Ker}}
\def\cu{{\rm cu}}
\def\c{{\rm c}}

\def\QQ{{\mathbb Q}}

\def\ZZ{{\mathbb Z}}

\def\GG{{\mathbb G}}

\def\DD{{\mathbb D}}
\def\FF{{\mathbb F}}

\def\cD{{\cal D}}
\def\cG{{\cal G}}
\def\cX{{\cal X}}
\def\cN{{\cal N}}
\def\cH{{\cal H}}
\def\cA{{\cal A}}
\def\cY{{\cal Y}}

\def\cI{{\cal I}}
\def\cC{{\cal C}}
\def\cM{{\cal M}}
\def\cZ{{\cal Z}}

\def\cV{{\cal V}}
\def\cW{{\cal W}}
\def\cQ{{\cal Q}}

\def\cO{{\cal O}}
\def\cU{{\cal U}}
\def\cF{{\cal F}}
\def\cP{{\cal P}}
\def\cL{{\cal L}}
\def\cK{{\cal K}}
\def\cJ{{\cal J}}
\def\cB{{\cal B}}

\def\B{{\hfill$\square$}}

\def\va{{\varphi}}

 \makeatletter
\renewcommand{\subsection}[1]{\@startsection{subsection}{2}{\z@}
            {-3.25ex plus -1ex minus -.2ex}{-1sp}{\normalsize\bf}
             {\ignorespaces#1 }}
 
\makeatother

\begin{document}
\sloppy

\title{Foliations in moduli spaces of abelian varieties}
\author{Frans Oort}
\date{Version 21-XI-2003}
\maketitle
\n
This paper appeared: J. Amer. Math. Soc. 17 (2004), 267-296.

\section*{Introduction}
\n
In this paper we study   abelian varieties and  $p$-divisible 
groups in characteristic $p$. 

Even though a non-trivial
deformation of an abelian variety can produce a non-trivial Galois-representation, 
say on the Tate-$\ell$-group of the generic fiber (in any characteristic), 
the geometric generic fiber has a ``constant'' Tate-$\ell$-group. 
The same phenomenon we encounter for the $p$-structure in positive characteristic: 
any two {\it ordinary} 
abelian varieties of the same dimension over an algebraically closed field have isomorphic
$p$-divisible groups.

However, for non-ordinary abelian varieties this seems to break down. The fascinating structure
which comes out of this is that the maximal locus where a given geometric isomorphism class of a
$p$-divisible group is realized (e.g. in a family of abelian varieties) is a locally closed set. The locus defined by the geometric isomorphism type of a $p$-divisible group
will be called a ``central leaf'', see \ref{cl}, and \ref{defcl}. This gives rise 
to a ``foliation'' of the open stratum attached to a Newton polygon $\xi$; the dimension of 
any ``leaf'' in the same Newton polygon stratum 
solely depends on $\xi$. In  extreme cases 

either the leaf is the whole stratum, 
as in the ordinary case, 

or in the ``almost ordinary case'' (the $p$-rank equals $g-1$),

or a leaf is zero-dimensional as in the supersingular case;\\
in intermediate cases a leaf can be a proper subset and still be positive dimensional: 
we have worked out the example of $g=4$ in \ref{exa}.

\vn
Between two central leaves in the same Newton polygon stratum 
there is a correspondence by iterated $\alpha_p$-isogenies. 
This leads naturally
to the study of ``Hecke-orbits'' by such isogenies. In the case of moduli spaces of abelian 
varieties we obtain ``isogeny leaves'', see \ref{isogl}. The beauty of these two structures is 
that they are 
almost transversal, that both feel like a foliation, and that they give, 
up to a finite morphism, 
a natural product structure on every irreducible component $W$ of an open Newton polygon stratum, 

\vn
We prove the following theorems:
\begin{description}
\item[\ref{fconst}]  -  a finite level of a \gfc \ family of  p-divisible 
groups becomes constant over an appropriate finite cover of the base;

\item[\ref{C}, \ref{Cp}] - a central leaf is closed in its open Newton polygon stratum;

\item[\ref{dim=}] - central leaves are smooth over the base field, and  the  dimension 
of a central leaf depends only on the Newton polygon;

\item[\ref{isogl}] - for polarized abelian varieties, the union of all irreducible iterated 
$\alpha_p$-Hecke-orbits through one point is a closed subset;  

\item[\ref{TI}] -  every component of a Newton polygon stratum is up to a finite morphism 
isomorphic  with the product of any of the isogeny leaves with a finite cover of any of the central leaves. 
\end{description}

\vn
{\bf Motivation / explanation.} In the moduli space of abelian varieties in characteristic 
$p$ we  consider Hecke orbits related to isogenies of degree prime to $p$, and 
Hecke orbits related to iterated $\alpha_p$-isogenies. The first ``moves'' points in 
a central leaf: under such isogenies  the geometric $p$-divisible group does not change;
the second moves points in an isogeny leaf. We can expect that these two natural foliations 
describe 
these two ``transversal'' actions: see \ref{HO} and \ref{HaO}.

\vn
{\bf EO-strata and central leaves.} In \cite{EO} we described the stratification of 
the moduli space of abelian varieties given by $p$-kernels: $(A,\lambda)$ and $(B,\mu)$
are in the same EO-stratum if $(A,\lambda)[p] \cong (B,\mu)[p]$, the isomorphism over an
algebraically closed field.  It is natural to consider more generally the relation given by 
the isomorphism type of $(A,\lambda)[p^i]$ for some $i \in \ZZ_{\geq 1}$; for every  
$i$ this defines
naturally  subsets of the moduli space, but we should not expect only finitely many of such 
subspaces. The central leaves studied in this paper
are the ones obtained by choosing $i$ large enough (depending on $g$), which is the same,
see \ref{2}, as 
considering the isomorphism type of $(A,\lambda)[p^{\infty}]$.

\vn
{\bf Isogeny correspondences.} Such correspondences are finite-to-finite in characteristic 
zero, or, more generally if only isogenies of degree prime to the characteristic are considered.
In positive characteristic such correspondences in general blow up and down.
{\it However it turns out that restricted to central leaves  all correspondences are 
finite-to-finite},
see \ref{icorr}.

\vn
{\bf Acknowledgments.} {\footnotesize I thank Thomas Zink for his ideas; 
suggestions by him have been incorporated in this paper; especially
his definition of \csd \ $p$-divisible groups, his paper \cite{Z2} 
and our cooperations on \cite{FO.Z} have been of value for me in 
preparing this paper. 
I thank Johan de Jong for discussions while preparing this paper, 
and for careful reading of the manuscript.  
Some aspects of this paper already have been incorporated by 
Elena Mantovan in  her Harvard PhD-thesis , see \cite{EM}, after discussions we had 
in the Fall of 2000 and in 
the Spring of 2002. I am grateful to MIT, where part of this paper was written. }  

\vn
Some terminology and notations used in this paper are brought together in Section \ref{Not} 
and Section \ref{T}. In particular a field denoted by $k$ will be supposed {\it algebraically closed}, 
$k = \overline{k} \supset \FF_p$.

\section{Preliminaries on $p$-divisible groups and on finite group schemes}
\subsection{}\label{defgfc}{\bf Definition.} {\it Let $S$ be a scheme, and let 
$X \to S$ be a $p$-divisible group. We say that $X/S$ is } \gfc \   
{\it if there exist a field $K$, 
a $p$-divisible group $X_0$ over $K$, a morphism $S \to \Spec(K)$, 
and for every $s \in S$  
an algebraically closed 
field $k \supset \kappa(s) \supset K$ containing  the residue class 
field of $s$  and an isomorphism $X_0 \otimes k \cong_k  X_s \otimes k$.}\\
The analogous terminology will be used for (polarized) abelian schemes.

\subsection{}\label{defconst}{\bf Definition.} {\it Let $T$ be a scheme and let 
$G \to T$ be a  group scheme. We say that $G$} is constant over $T$ 
{\it if there exist a field $L$, a group scheme $G_0$ over $L$, a morphism $T \to \Spec(L)$ 
 and an isomorphism} 
$$G_0 \times_{\Spec(L)} T \quad\cong_T\quad G.$$

\vn
In this section we show:
\subsection{}\label{fconst}{\bf Theorem.} {\it Let $S$ be a   scheme which satisfies condition} (N), see \ref{N}; {\it   
let $\cX \to S$ be a $p$-divisible group; let $n \in \ZZ_{\geq 0}$. Suppose that
$\cX \to S$ is} \gfc. {\it  Then there exists a finite surjective morphism
$T_n = T \to S$, such that $\cX[p^n] \times_S T$ is constant over $T$.}

\vn
{\bf Remark.} The choice of $T \to S$ can be made in such a way that this equals a composition of finite morphisms $T \to T' \to S' \to S$, where $S' \to S$ factors the normalization, $\tilde{S} \to S' \to S$, and  $T' \to S'$ is purely inseparable and $T \to T'$ is  \'etale.

\vn
For the proof of \ref{fconst} we need some preliminaries.\\ 
{\it We say that the condition $\ast(Z/S,n)$ is  satisfied, where   
$Z \to S$ is a $p$-divisible group, and $n \in \ZZ_{\geq 0}$, 
if there exists a finite surjective morphism $T \to S$ such that 
$Z[p^n] \times_S T$ is constant.} 

\subsection{}\label{1}{\bf Lemma.}  {\it For any completely 
slope divisible $Y \to S$, where $S$ is noetherian and integral,   the property  $\ast(Z/S,N)$ is satisfied for every $N \in \ZZ_{\geq 0}$.}\\
{\bf Proof.} At first we claim that {\it there exist a purely inseparable, finite, 
surjective morphism 
$T_1=T \to S$ such that} $Y[p^N]_T = \oplus_i ((Y_i/Y_{i-1})[p^N])_T$, 
in  the notation of \ref{defcsd}.

In fact,  for large $r \in \ZZ$,  over $T' = S^{(p^{-rs})}$ the homomorphism 
$$\left(\psi:= \frac{F^s}{p^{s-t_i}}\right)^r: 
((Y_i/Y_{i-1})[p^N])^{(p^{-rs})}_{T'} \longrightarrow Y_i[p^N]_{T'} $$
splits off this subquotient, because $\psi$ is an isomorphism on $Y_i/Y_{i-1}$ and nilpotent on $Y_{i-1}[p^N]$. Note  that $S$ is integral and noetherian, and $Y_i[p^N] \to S$ is finite; hence there exists a factorization 
$T' = S^{(p^{-rs})} \to T_1 \to S$, where  $T_1 \to S$ is finite, and for which  the claim holds. 

By \cite{FO.Z}, 1.10 we see that all $G_i =  (Y_i/Y_{i-1)}[p^N]$ become constant, for $1 \leq i \leq m$,
under a (separable) finite, surjective morphism $T_2 \to S$. Any $T \to S$
dominating $T_1 \to S$ and $T_2 \to S$ has the desired property. 
This proves the lemma.      \B 

\vn
For $p$-divisible groups $X, Y$  over a field and for $N \geq n \geq 0$ we denote by
$$\Phi_n: \Hom(X,Y) \to  \Hom(X[n],Y[n])  \quad\mbox{and}\quad \Phi^N_n:  
\Hom(X[N],Y[N]) \to \Hom(X[n],Y[n])$$
the natural restriction maps.
 
\subsection{}\label{bound}{\bf Lemma.} {\it For $p$-divisible groups $X$ and $Y$ over $k$, 
there exists for every $n \in \ZZ_{\geq 0}$ an integer $N(X,Y,n)$ such that for
every $N \geq  N(X,Y,n)$ we have  
$$\Ima(\Phi_n) = \Ima(\Phi^N_n).$$}
{\bf Proof.} For every $m \in \ZZ_{\geq 0}$ consider the functor 
$T \mapsto \Hom(X[p^m]_T,Y[p^m]_T)$ on $k$-schemes. As is easily seen
this functor is representable by a group scheme of finite type over $k$,
notation:
$G_m = \underline{\Hom}(X[p^m],Y[p^m])$. For $m \geq n$ the natural restriction
map yields a homomorphism, of algebraic groups 
$\rho^m_n: G_m \to G_n$.
We write 
$G_m' = (G_m)_{\rm red}$.  Note that the image of a homomorphism between
algebraic groups is a (closed) subgroup; consider the finite set
$\Phi_n(\Hom(X,Y)) \subset G'_n(k)$ as an algebraic subgroup of $G'_n$.
We obtain 
a descending chain  of algebraic groups 
$$G'_n \supset \rho^{n+1}_n(G'_{n+1}) \supset \cdots \supset  
\rho^{n+i}_n(G'_{n+i}) \supset \cdots \supset \Phi_n(\Hom(X,Y)).$$
The sequence of subgroups schemes $\{\rho^{n+i}_n(G'_{n+i}) \mid i \in \ZZ_{\geq 0}\}$
stabilizes on the one hand, and it is equal to $\Phi_n(\Hom(X,Y))$ on the other hand. \B

\vn
We are going to show that  $n \mapsto N(X,Y,n)$ can be chosen uniformly 
 depending only on the heights of the $p$-divisible groups in consideration.

\subsection{}\label{Nn}{\bf Proposition.} {\it For every $h \in \ZZ_{>0}$ there exists 
an integer $N(h',h'',n)$
such that for $p$-divisible groups $X$ and $Y$ over an algebraically closed 
field $k$ of height} $h'=$height$(X)$ {\it respectively} $h''=$height$(Y)$ 
{\it and for every $N \geq N(h',h'',n)$
the following images are equal:}
$$\Ima\left(\Phi_n\right)  = \Ima\left(\Phi^N_n\right).$$
{\bf Proof.} We split up the proof in various steps.\\
{\bf (1)} {\it Given $h$ there exists an integer $d(h)$ such that for every 
$p$-divisible group $X$ of height $h$ and with Newton polygon $\cN(X) = \beta$ there is 
an isogeny $\rho: H(\beta) \to X$ of degree $\deg(\rho) = d(h)$.}\\
This follows from \cite{Manin}, (3.4) and (3.5) on page 44, and from the property
that $H(\beta)$ admits endomorphisms of degree $p$.  \B(1)

\vn
{\bf (2)} In order to prove the proposition is suffices to show:\\
{\it For every $h \geq 0$ there exists $N(h)$ such that for every $p$-divisible group $Z$
of height $h$ and for every $N \geq N(h)$ we have: }
$$\Ima(\Phi_n: \End(Z) \to \End(Z[n])) = \Ima(\Phi^N_n: \End(Z(N)) \to \End(Z[n])).$$ 
Indeed, apply this result for $Z = X \oplus Y$.  \B(2)

\vn
{\bf (3)} {\it Let $H$ and $Z$ be $p$-divisible groups and let $\rho: H \to Z$ be an
isogeny of degree $\deg(\rho) = p^s$. Suppose that the property mentioned in} (2) 
{\it holds for $H$ with the function $b \mapsto N_H(b)$; then that property holds for
$Z$ with the function $n \mapsto N_Z(n) := N_H(n+s) + s$.} \\
{\bf Proof of (3).} We consider Dieudonn\'e modules and induced maps:
$$\DD(\rho: H \to Z) = \left(\DD(H) =: Q \hookrightarrow M := \DD(Z) \right).$$
For $n \geq 0$ and $N \geq N_H(n+s) + s$ we consider the inclusions:
$$p^NM \quad\subset\quad  p^{N_H(n+s)}Q\quad\subset\quad  p^{n+s}Q \quad\subset\quad   
p^{n+s}M \quad\subset\quad  p^nQ \quad\subset\quad Q  
\subset M \subset p^{-s}Q .$$
Suppose $\va_N \in \End(Z[N])$ restricts to $\va_n = \Phi^N_n(\va_N)$; we claim 
that in this case $\va_n$ can be lifted to $\va \in \End(Z)$. Indeed, consider 
$\psi_{n+s} \in \End(Z[p^{n+s}]$ defined by $\psi_{n+s} = p^s\va_n$. We denote
the endomorphisms induced on the Dieudonn\'e modules by the same symbols;
note that $\psi_{n+s}: M/p^{n+s}M \to M/p^{n+s}M$
restricts to $\psi'_{n+s}: Q/p^{n+s}M \to Q/p^{n+s}M$ by:
$$\frac{Q}{p^{n+s}M}  \quad\hookrightarrow\quad  \frac{M}{p^{n+s}M} 
\quad\stackrel{p^s\va}{\longrightarrow}\quad  
\frac{p^sM}{p^{n+s}M}   \quad\hookrightarrow\quad \frac{Q}{p^{n+s}M}.$$
We see that $\va_N$ restricts to 
$\psi'_{N_H(n+s)} \in \End(H[p^{N_H(n+s)}]) = \End(Q/p^{N_H(n+s)}Q)$;
hence the restriction to $\End(H[p^{n+s}])$ can be lifted to $H$. This shows 
that $\psi'_{n+s}  \in \End(Q/p^{n+s}M)$ can be lifted to $Q$. {}From this it follows 
that $\va_n \in \End(M)$ can be lifted to $\End(M)$.  \B(3)

\vn
{\bf (4)} {\it For a given $H = H(\beta)$ there is a function $N_H$ such that} 
(2) {\it is satisfied for $H$ and this function.
Hence, for all minimal groups of height $h$ there is a function 
$n \mapsto N_{\rm min}(h,n)$ satisfying } (2) {\it for all minimal groups of height $h$.}\\
The first statement is a consequence of \ref{bound}. The second follows because
there are only finitely many Newton polygons of a given height.  \B(4)  

\vn
{\bf Remark} (not used in this paper). {\it An easy and explicit consideration shows that  
for any Newton Polygon $\beta$
the function $N_{H(\beta)}(n) := n+1$ satisfies the property in} (4).

\vn
We have seen in (2) that it suffices to show the proposition for endomorphisms. 
For a given $h \in \ZZ_{>0}$ we define $N(h,n) = N_{\rm min}(h,n+s) + s$, with  
 with $s := d(h)$ as in (1).
By (1) we see that for every $Z$ of height $h$ there exists an isogeny 
$H = H(\cN(Z)) \to Z$ of degree $d(h)$. By (4) we know that property (2) holds for the
function $N_H(m) = N_{\rm min}(h,m)$. By (3) we see that property (2) holds for $Z$ with the function
$N_Z(n) = N_{\rm min}(h,n+s) + s$.   \B\ref{Nn}  

\vn
Note that the detour via minimal groups is not necessary: the proposition follows,
using special Dieudonn\'e submodules as in \cite{Manin},  
with a less explicit bound, from:
\cite{Manin}, (3.4) and (3.5) on page 44, \ref{bound}  and (2) and (3). Note
however that the number of \csd \ $p$-divisible groups, up to isomorphism over $k$,
of given height, in general is not finite; for that reason we dit not use these
$p$-divisible groups in the previous proof.

\vn
The following result will be used several times in this paper. Probably this
fact is known to all experts in the field.

\subsection{}\label{2}{\bf Corollary.} {\it For any $h \in\ZZ_{\geq 0}$
there exists $N(h) \in \ZZ_{\geq 0}$ with the property that for 
$p$-divisible groups $X_1$, and $X_2$ of height $h$ over  
an algebraically 
closed field, and any $n \geq N(h)$, then:
$$X_1[p^n] \quad\cong\quad X_2[p^n] \quad\Rightarrow 
\quad X_1 \quad\cong\quad X_2.$$ }
Apply \ref{Nn}  with $n=1$.  \B 

\subsection{}\label{3}{\bf Corollary.} {\it  Any completely 
slope divisible $Y \to S$, over   a noetherian and integral $K$- scheme $S$ 
with $S(K) \not= \emptyset$ is
\gfc.}\\
This follows from \ref{1} and \ref{2}.  \B

\subsection{}\label{4}{\bf Lemma.} {\it Suppose given $X \to S$, 
a $p$-divisible group, and $i \in  \ZZ_{\geq 0}$, the functor
$$T \mapsto Gr(T) := \{G \hookrightarrow X_T\mid \rk(G/T) = p^i\}$$
is representable by a proper scheme $Gr_{X/S,i} = Gr \to S$,
and a universal $\cG \to Gr$.}\\
{\bf Proof.} There is a locally free sheaf $\cB$ of algebras on $S$ giving 
$X[p^i] \to S$. Any $G \hookrightarrow X_T$ as considered is included
in $X[p^i]_T$, and it is given by an ideal sheaf $\cI$ in $\cB \otimes \cO_T$.
Hence the functor is a subfunctor of the related Grassmannian of $\cB/S$. The
properties that $\cB \to \cB/I$ with $G = \Spec(\cB/\cI)$ is in fact a
morphism of bialgebras is a closed condition on  points in that Grassmannian.  \B 

\subsection{}\label{5}{\bf Lemma.} {\it Let $Z_1$ and $Z_2$ be $p$-divisible 
groups over an algebraically closed field $k$, and let $i \in \ZZ_{\geq 0}$.
We write
$$\Psi = \Psi(Z_1,Z_2; i) := \{G \subset Z_1\mid \rk(G) = p^i, \ \ \exists \psi: Z_1 \to Z_2,
\ \  \Ker(\psi) = G
\};$$
then $$\#\left(\Psi\right) \quad<\quad \infty.$$  }\\
{\bf Proof.} The module $\Gamma := \Hom(Z_1,Z_2)$ is free of finite rank over $\ZZ_p$.
Hence $\Gamma/p^i\Gamma$ is finite. Let 
$\Phi' \subset \Gamma$ and its image $\Phi''\subset \Gamma/p^i\Gamma$
    be defined as:
$$\Phi' := \{\psi \in \Gamma\mid \rk(\Ker(\psi)) = p^i\}; \ \  \Phi' \ \mbox{\rm mod}\ \ p^i\Gamma   =: \Phi''\subset \Gamma/p^i\Gamma.$$
If $\psi_1 \in \Phi'$, and $\rho \in \Phi'$, and $\psi_1 + p^i{\cdot}\rho =: \psi_2 \in \Phi'$, then $\Ker(\psi_1) = \Ker(\psi_2)$.
The map $\Phi' \twoheadrightarrow \Phi$  defined by   $\psi \mapsto  \Ker(\psi)$ induces a surjective map 
$\Phi'' \twoheadrightarrow \Psi$; hence  $\#\left(\Psi\right) < \infty.$ 
\B

\subsection{}{\bf Proof of Theorem \ref{fconst}.}  Let $S' \to S$ be 
the normalization 
morphism, and consider $\cX' = \cX_{S'} \to S'$; by \cite{FO.Z}, Theorem (2.1),
or \cite{FO.Z}, Proposition (2.7), 
we know there 
exists an isogeny $\va: \cY \to \cX'$ over $S'$, where $\cY/S'$ is \csd.  
Let $p^i := \deg(\va)$. 
We are going to show $\ast(\cX_{S'},n)$. We choose  $N \in \ZZ_{\geq 0}$ 
such that 
$N \geq n+i$. As $\cY \to S'$ satisfies
$\ast(\cY,N)$, see \ref{1}, we can choose a finite $T \to S'$, with $T$ connected,
a field $L$,   a morphism $T \to \Spec(L)$,
a $p$-divisible 
group $Y_0$ over  $L$ , 
 and an isomorphism
$$Y_0[p^N] \times T \quad\stackrel{\sim}{\longrightarrow}\quad \cY[p^N] \times_{S'} T;$$
from now on we fix this isomorphism and write it as an equality.
We write $X_0 = (\cX_{T})_0$, and we have the morphism $\va_{T,0}: Y_0 \to X_0$. 
Because $N \geq i$ we see that $\Ker(\va) \subset \cY[p^i] \subset \cY[p^N]$; 
hence $Gr := Gr_{\cY,i} =  Gr_{Y_0,i} \times T$. Consider the finite set 
$\Psi = \Psi(Z_1 = Y_0, Z_2 = X_0; i)$ as in \ref{5}.

The morphism $\va_T: \cY_{T} \to  \cX_{T}$  defines $\cG := \Ker(\va_{T}) 
\subset \cY[p^N]_T = Y_0 \times T$; 
this is equivalent with a section 
$$Gr \quad\stackrel{\stackrel{\sigma}{\longleftarrow}}{\longrightarrow}\quad T.$$ 
Note that $\cY_T/T$ is \csd; hence for every geometric point $t$ of $T$ 
there exists some isomorphism $\cY_t \cong Y_0 \otimes k$, see \ref{3}; 
also we have some isomorphism $\cX_t \cong X_0 \otimes k$. 
Hence for every geometric 
point $t$ of $T$ we have $\sigma(t) \in \Psi = \Psi(Y_0,X_0;i)$; note that  
this statement does not
depend on the choice of the isomorphisms $Y_t \cong Y_0 \otimes k$
and $X_t \cong X_0 \otimes k$: note that $\Aut(Y_0)$
acts from the right on $\Psi$, and $\Aut(X_0)$ acts from the left on $\Psi$. 
As $T$ is connected, and $\sigma(t) \in \Psi$ for every $t$, this shows 
$\sigma: T \to Gr_{Y_0,i} \times T$ to be
a constant map; hence there exists $G_0 \subset Y_0$ such that  
$$\Ker((\va_T)_0) = G_0 \times T \subset \cY[p^N]_T = Y_0 \times T.$$
As $N \geq n+i$ we see that  
$$\va_T(\cY[p^N]_T) \supset \cX[p^n]_T, \quad\mbox{hence}\quad 
\cY[p^N]_T \supset \cH:= \va_T^{-1}(\cX[p^n]_T) = \Ker(p^n{\cdot}\va_T).$$
We have 
$$\cX[p^n] = \cH/\cG, \quad\mbox{and}\quad G_0 \subset H_0 \subset Y_0.$$
We see that $\cH \subset \cY[p^N]_T = Y_0 \times T$ defines a section $\Sigma$ of 
$Gr_{Y_0,j} \times T$, with $j=hn+i$. Every geometric fiber $\cH_t$ 
is the kernel of $p^n{\cdot}\va_t:  Y_0 \cong \cY_t \to \cX_t \cong X_0$; hence $\cH$ 
defines 
a section $\Sigma: T \to Gr= Gr_{Y_0,j} \times T$; this section has  
the property that for every geometric point t we have 
$\Sigma(t) \in \Psi(Y_0, X_0; j)$;
hence, by the same argument as before,
$$\cH = H_0 \times T \subset \cY[p^N]_T = Y_0 \times T, \quad\mbox{with}\quad 
\cG \subset \cH \quad\mbox{giving}\quad  (\cG \subset \cH) = (G_0 \subset H_0) \times T.$$     
This proves that
$$\cX[p^n]_T = \cH/\cG = (H_0/G_0) \times T.$$ 
This proves $\ast(X,n)$ and we have concluded the proof of the theorem.
\hfill$\Box$\ref{fconst}

\section{Central leaves}
\subsection{}\label{cl}{\bf Notation.} {\it Let $K \supset \FF_p$ be a field, and let $X$ be a 
$p$-divisible group over $K$. Let $S \to \Spec(K)$ be a noetherian scheme over $K$, 
and let $\cY \to S$ 
be a $p$-divisible group over $S$. We write:}
$$\cC_X(S) \quad:=\quad  \{s \in S \mid \exists k = \overline{k} \supset \kappa(s), \ \ 
\exists \cong: \ \cY_s \otimes_{\kappa(s)} k \cong X \otimes _K k\}.$$
Note that $\cY_s \otimes_{\kappa(s)} k \cong X \otimes _K k$ as in the definition exists iff there is a field extension $\kappa(s) \subset L$  and an isomorphism $\cY_s \otimes_{\kappa(s)} L \cong X \otimes _K L$.\\
{\bf Remark.}  We will write $\cC_X(\cY \to S)$ if confusion might be possible. The isomorphism in the definition exists iff such an isomorphisms exists over an algebraic closure of $\kappa(s)$. We have defined 
$\cC_X(S) \subset S$ as a {\it subset}. Once  \ref{C} has been proved, working over 
a perfect field, we consider $\cC_X(S) \subset S$ as a closed subscheme of $\cW^0_{\cN(X)}(S)$ 
with the induced reduced scheme structure. Once the theorem has been proved, in such situations
the formation $\cC_X(-)$ commutes with base change.

\vn 
We remind the reader of the following notation. Let $\beta$ be a Newton polygon, see \ref{NP}, and let 
$\cY \to S$ be a $p$-divisible group with the same $d$, and $h = d+c$. We write:
$$\cW_{\beta}^0(S) = \{s \mid \cN(\cY_s) = \beta\},$$
and
$$\cW_{\beta}(S) = \{s \mid \cN(\cY_s) \prec \beta\}.$$ 
By Grothendieck-Katz, see \cite{Katz}, 2.3.1 and 2.3.2, we know that 
$\cW_{\beta}(S) \subset S$ is 
closed and $\cW_{\beta}^0(S) \subset S$ is locally closed; these sets 
are given the induced, reduced scheme structure. Note that  
$\cC_X(S) \subset \cW_{\cN(X)}^0(S)$.

\vn
In this section we will show:
\subsection{}\label{C}{\bf Theorem.} {\it Let $K$ be a field, let $X$ be a 
$p$-divisible group over $K$. Let $S \to \Spec(K)$ where $S$ satiefies is an excellent scheme over a field $K$. 
Let $\cY \to S$ be a $p$-divisible group. Then
$$\cC_X(S) \subset \cW_{\cN(X)}^0(S)$$
is a} closed {\it subset.}

\subsection{}{\bf Definition.}  {\it  Let $S$ be a scheme. Let $G_1 \to S$ and $G_2 \to S$ 
be finite, locally free group schemes over $S$. Define the ``isom - functor'' $\cI$
by:}
$$\forall T \to S, \ \ \  \cI(T) = \Isom(G_1 \times_S T, G_2 \times_S T).$$

\subsection{}\label{Isom}{\bf Lemma.} {\it This functor is representable by 
a scheme $I(G_1, G_2) = I \to S$
of finite type over $S$.}\\
{\bf Proof.} Let $\cB_1$, respectively $\cB_2$ be the $\cO_S$-algebra  defining
$G_1 \to S$, respectively $G_2 \to S$. An isomorphism as considered defines
an $\cO_T$-isomorphism $\cB_2 \otimes \cO_T \to \cB_1 \otimes \cO_T$; the functor of 
such maps is represented by a scheme of finite type over $S$. The condition that 
such maps moreover define an isomorphism of bialgebras (defining an isomorphism of 
group schemes) is a closed condition in this scheme.   \B

\subsection{}\label{U}{\bf Corollary.} {\it With the notation as in} \ref{C}, 
{\it by} \ref{Isom} {\it we see that
$\cC_X(S) \subset S$ is a constructible subset. Hence it contains a dense  
subset $U \subset \cC_X(S)$, open in the closure of $\cC_X(S)$.}\\
Indeed, use \ref{2}.  \B

\subsection{}{\bf Proof of \ref{C}.} Let 
$S' = \mbox{Zar}\left(\cC_X(S) \subset \cW^0_{\cN(X)(S)}\right)$
be the Zariski-closure. Let $S'' \to S'$ be the normalization, and 
$\cX'' = \cX \times _S S''$. Note that $S'' = \cW^0_{\cN(X)}(S'')$: the Newton
polygon of $\cX'' \to S''$ is constant.
By \cite{FO.Z}, Theorem 2.1 there exists 
$\va: \cY \to \cX''$ such that $\cY \to S''$ is \csd. Let  $\deg(\va) = p^i$. Note that 
$\cY \to S''$ is \gfc, see \ref{3}. Consider $N(h)$ as in \ref{2}, where $h$ is 
the height of $\cX/S$. Choose $N \in \ZZ_{>0}$  such that $N \geq N(h) + i$.
As $\cY \to S''$ satisfies $\ast(\cY,N)$, see \ref{1}, we can choose a finite surjective 
$T \twoheadrightarrow S''$, a field $L$, a morphism $T \to \Spec(L)$, a point $0 \in T(L)$,
and an isomorphism
$$Y_0[p^N] \times T \quad\stackrel{\sim}{\longrightarrow}\quad \cY[p^N] 
\times_{S''} T;$$
we fix this isomorphism, and write it as an equality.

The morphism  $\va_T: Y_T \to X_T$ defines $\cG := \Ker(\va_{T}) 
\subset \cY[p^N]_T \subset Y_0 \times T$, and this is equivalent with a section
$$Gr_{\cY,i} = Gr_{Y_0,i} \times T \quad\stackrel{\stackrel{\sigma}{\longleftarrow}}
{\longrightarrow}\quad T.$$
Choose $U \subset \cC_X(S)$ as in \ref{U}; as $\cX \times_{S} (U \times_S T)$
and $\cY \times_{S''} (U \times_S T)$ are \gfc \ we conclude that
$\sigma$ is a constant map over $U \times_S T$; hence it is constant over $T$,
hence there exists $G_0 \subset Y_0$ such that
$$\Ker((\va_T)_0) = G_0 \times T \subset \cY[p^N]_T = Y_0 \times T.$$ 
As $N \geq N(h) + i$,  defining  $\cH := \va_T^{-1}(\cX[p^{N(h)}_T])$,
we see, as in the proof of \ref{fconst}, that 
$$\mbox{for every geometric point}\quad \  \  t \in T(k) \quad\mbox{we have}\quad \ \ 
\cX_t[p^{N(h)}] \cong X_0[p^{N(h)}] \otimes k.$$
By \ref{2} this proves that $\cX_t \cong X_0 \otimes k$ for every geometric
point $t \in T(k)$. Hence $\cX_T \to T$ is \gfc; this proves that $\cC_X(S) = S'$,
which finishes the proof of Theorem \ref{C}.
\hfill$\Box$

\subsection{}\label{cu}{\bf Proposition / Notation} $\cu(-)$. {\bf (i)} 
{\it Let $\va_0: Y_0 \to X_0$ be an isogeny of $p$-divisible groups over a field $k$; 
let $R$ be a complete local noetherian domain with residue class field $k$  
and let $Y \to T := \Spec(R)$
be a \gfc \ $p$-divisible group. Consider $C:= \cC_{X_0}(\Def(X_0))$ and the restriction
$X \to C$ of the universal family.
There exists a correspondence 
$T \leftarrow \Gamma \to  C = \cC_X(\Def(X))$ and an isogeny $\va: Y_{\Gamma} \to X_{\Gamma}$
extending $\va_0$ with $T \leftarrow \Gamma$ a finite, surjective morphism.}\\
{\bf (ii)} {\it For every Newton polygon $\beta$ there exists 
an integer $\cu(\beta) \in \ZZ_{\geq 0}$ } (c = central leaf, u = unpolarized)
{\it such that for every $p$-divisible group $X$ with $\cN(X) = \beta$ the closed locus 
$\cC_X(\Def(X))$ is pure of dimension $\cu(\beta)$.}\\
The first part follows from \ref{fconst}; for more details see the proof of \ref{LemT}. The second part is an immediate 
consequence of the first.   \B

\subsection{}{\bf Remark.} In \cite{C.FO} we prove the following formula. Let $\beta$ 
be a Newton polygon of dimension $d$ and height $h$, given by the slopes 
$0 \leq \lambda_1 \leq \cdots \leq \lambda_h \leq 1$. Let $\beta^{\ast}$ be the (upper convex)
polygon given by $\lambda_h, \cdots , \lambda_1$. We have 
$\cu(\beta) \quad=\quad \sum_{0 < j < h}  \   \left( \beta^{\ast}(j) - \beta(j) \right).$

\section{Central leaves: the polarized case}
\subsection{} Let $A \to S$ be an abelian scheme. Let $\lambda: A \to A^t$
be a polarization on $A \to S$; this  isogeny 
is symmetric, i.e.  $\lambda = \lambda^t$ under the identification
$A \stackrel{\sim}{\longrightarrow} A^{tt}$.

Let $Y \to S$ be a $p$-divisible group. Let $Y^t \to S$ be its Serre-dual.
We say that an isogeny $\lambda: Y \to Y^t$ is a quasi-polarization on $Y/S$ if
$$(\lambda^t : Y^{tt} \to Y^t) \quad=\quad -(Y^{tt} = Y 
\stackrel{\lambda}{\longrightarrow} Y^t).$$
Note that a polarization (a symmetric isogeny) on an abelian scheme $A$ 
induces a quasi-polarization (an anti-symmetric isogeny) on its
$p$-divisible group $X = A[p^{\infty}]$, see \cite{Oda}, Proposition 1.12; 
a polarization on $A$ and the induced 
quasi-polarization on $X = A[p^{\infty}]$
will be denoted by the same symbol, unless there might be danger of confusion.

\subsection{}\label{pcl}{\bf Notation.} {\it Let $K \supset \FF_p$ be a field, and let $(X,\lambda)$ be a 
quasi-polarized $p$-divisible group over $K$. Let $S \to \Spec(K)$ be a scheme over $K$, 
and let $(\cY,\mu) \to S$ 
be a quasi-polarized $p$-divisible group over $S$. We write:}
$$\cC_{(X,\lambda)}(S) \quad:=\quad  \{s \in S \mid \exists k = \overline{k} \supset \kappa(s),
 \ \ 
\exists \cong: \ (\cY_s,\mu_s) \otimes_{\kappa(s)} k \cong (X,\lambda) \otimes _K k\}.$$
In case $K$ is perfect we consider $\cC_{(X,\lambda)}(S)$ as a scheme with the induced reduced scheme structure.

\subsection{}\label{Cp}{\bf Theorem.} {\it Let $K$ be a field, let $(X,\lambda)$ be a 
quasi-polarized $p$-divisible group over $K$. Let $S \to \Spec(K)$ where $S$ is an excelent scheme over $K$. 
Let $(\cY,\mu) \to S$ be a quasi-polarized $p$-divisible group over $S$. Then
$$\cC_{(X,\lambda)}(S) \subset \cW_{\cN(X)}^0(S)$$
is a} closed {\it subset, which is a union of connected components of $\cC_X(S)$.}\\
By \ref{C} we know that $\cC_X(S) \subset \cW_{\cN(X)}^0(S)$ is closed. Take the union of those irreducible components $C \subset  \cC_X(S)$ for which for the generic point $\eta \in C$ an isomorphism 
$(\cY_{\eta},\mu_{\eta}) \otimes_{\kappa(\eta)} k \cong (X,\lambda) \otimes _K k$ exists; using Theorem \ref{fconst} and Corollary \ref{2} we see that this union equals $\cC_{(X,\lambda)}(S)$.
\B

\subsection{}\label{defcl}{\bf Central leaves.} 
For any point $[(A,\lambda)] =x \in \cA_{g,d} \otimes \FF_p$  we choose a field $K$
over which $(X:=A[p^{\infty}], \lambda)$ is defined, and we consider:
$$C(x) := \cC_{(X,\lambda)}(\cA_{g,d} \otimes K).$$
Suppose $K$ is perfect; a geometrically irreducible component of $C(x)$
with the induced reduced scheme structure  will be called 
a {\bf central leaf}. We will see that for $x \in \cA_{g,d,n} \otimes k$ with $n \geq 3$
prime to $p$ there is exactly one irreducible component of $C(x)$ containing $x$; this component
will be denoted by $C_x \subset \cA_{g,d,n} \otimes k$.  See \ref{ConjCx}.

\subsection{}\label{ss}{\bf Notation.} Over $\FF_{p^2}$ we consider two series of quasi-polarized 
superspecial $p$-divisible groups (see \cite{Li.FO}, 6.1):

\vn
{\bf I}$_r$: $(S = G_{1,1},\tau_r)$, where on $M = \DD(S) = Wx \oplus W{\cdot}\cF x$, 
with $\cF x = \cV x$ the pairing $\tau_r$
is given by: 
$$\epsilon \in W - pW, \ \epsilon^{\sigma} = -\epsilon, \ 
\ <x,\cF x> = p^r{\cdot}\epsilon, \ \  \deg(\tau_r) = p^{2r} \ \ r \in \ZZ_{\geq 0};$$
{\bf II}$_r$: $(T = ((G_{1,1})^2,\nu_r)$ , where on 
$M = \DD(T) = Wx \oplus W{\cdot}\cF x \oplus Wy \oplus W{\cdot}\cF y$ 
the pairing is given by $<x,y>=p^r,  \ \ <\cF x,\cF y> = p^{r+1}$
(and all other pairs of base vectors pair to zero); 
here $\deg(\nu_r) = p^{4r+2}$, and $r \in \ZZ_{\geq 0}$.

\vn
In \cite{Li.FO}, 6.1 it is proved that over an algebraically closed field a superspecial
quasi-polarized $p$-divisible group is isomorphic to a direct sum of summands 
each of type I or II.

\subsection{}{\bf Notation.} Choose coprime, non-negative integers  $m, n \in \ZZ_{\geq 0}$,
with $m > n$. For any $r \in \ZZ_{\geq 0}$ we consider over $\FF_p$ a quasi-polarized 
$p$-divisible group
$(U,\zeta)$,  sometimes written as $U = U(m,n), \ \  \zeta = \zeta_r(m,n)$ by:
$U(m,n) = H_{m,n} \oplus H_{n,m},$ and
$$\zeta = \left(
\begin{array}{cc}
0 &  -\overline{\pi}^r\\
\pi^r & 0 
\end{array} 
\right) : \quad H_{m,n} \oplus H_{n,m}\quad\longrightarrow\quad  (H_{m,n} \oplus H_{n,m})^t,$$
where $\pi \in \End(H_{m,n} \otimes \FF_p)$ is a uniformizer in 
$\End(H_{m,n} \otimes \overline{\FF_p})$, see \cite{Purity}, 5.3,   
and  $\pi = \pi_{m,n}: H_{m,n} \to H_{m,n} = (H_{n,m})^t$, and
$\overline{\pi}=\pi_{n,m} \in \End(H_{n,m} \otimes \FF_p)$ is a uniformizer.
Note that $\pi^t = \overline{\pi}$ under the identification $(H_{m,n})^t = H_{n,m}$.

We assume gcd$(m_i,n_i) = 1$, and $m_i > n_i$,  and $n_i/(m_i+n_i) < n_j/(m_j+n_j)$  for $j>i$.

\subsection{}\label{eld}{\bf Proposition} (elementary divisors). {\it In this 
proposition we work over 
an algebraically closed field $k$. Suppose $\xi$ is a symmetric Newton polygon, 
and let $(H,\zeta)$ be a quasi-polarized $p$-divisible group with 
$$H \quad\cong \quad H(\xi) \quad=\quad (\oplus_i \ \  (H_{m_i,n_i})^{r_i}) 
\bigoplus (G_{1,1})^s 
\bigoplus(\oplus_i \ \  (H_{n_i,m_i})^{r_i}),$$} 
i.e. $H$ is minimal, see \ref{min}. {\it Then there exists an isomorphism
$$(H,\zeta)\quad \cong \quad 
\oplus_i\oplus_{1 \leq j \leq r_i}\left(U_{m_i,n_i},\zeta_{d_{i,j}}(m_i,n_i)\right) 
\bigoplus (\oplus_b(Q_b,\beta_b)),$$
where the non-negative integers $d_{i,1} \mid \cdots \mid d_{i,r_i}$ 
are  the elementary divisors of the 
quasi-polarization  on this isoclinic part, and  every $(Q_b,\beta_b)$ is of type {\bf (I)}
or of type {\bf (II)} as in} \ref{ss}. \\
{\bf Proof.}  It suffices to show this for the supersingular part, and for the parts
$(H_{m_i,n_i}\oplus H_{n_i,m_i})^{r_i}$ and $(G_{1,1})^s$ separately: $\zeta$ is in block form. 
For the supersingular part this is proven in \cite{Li.FO}, 6.1, proposition. 
A quasi-polarization  $\zeta$ on $(H_{m_i,n_i} \oplus H_{n_i,m_i})^{r_i}$ is in block form
$$ \zeta = \left(
\begin{array}{cc}
0 &  \gamma\\
\beta & 0 
\end{array}\right), \quad \beta: H_{m_i,n_i}^{r_i} \to ((H_{n_i,m_i})^t)^{r_i}, 
\quad \gamma: (H_{n_i,m_i})^{r_i} \to (H_{m_i,n_i}^t)^{r_i};$$
by $\zeta^t = -\zeta$, on $p$-divisible groups, we conclude that this is a quasi-polarization iff $\beta^t = -\gamma$.
We use 
that $E = \End(H_{m,n} \otimes k)$ is the maximal order in  $\End^0(H_{m,n} \otimes k)$,
and the same for $H_{n,m}$: by choosing an appropriate basis for $(H_{n_i,m_i})^{r_i}$ 
and for $(H_{m_i,n_i}^t)^{r_i}$ the morphism $\beta \in {\rm M}(r_i,E)$ can be diagonalized 
with only powers of $\pi$ on the diagonal, see \cite{Reiner}, Th. 17.7.  This finishes the
proof of the proposition.
\hfill$\Box$

\subsection{}\label{unique}{\bf Corollary.} {\it Over an algebraically closed field 
there is, up to isomorphism, precisely one principal quasi-polarization on $H(\xi)$.}
\hfill$\Box$

\subsection{}{\bf Corollary.} {\it Let $X$ be a $p$-divisible group over $k$. The number
of principal quasi-polarizations on $X$ up to isomorphism is finite.}\\
{\bf Proof.} Choose an isogeny $\psi: H = H(\cN(X)) \to X$; write $\deg(\psi) = p^i$. 
Any principal polarization  $\lambda$
on $X$ pulls back to a polarization $\psi^{\ast}(\lambda)$ of degree $p^{2i}$. 
The number of  polarizations of a given degree on $H$  is finite by \ref{eld}.
The isogeny $\psi$ gives an identification $\End^0(X) \to \End^0(H)$ by $\gamma \mapsto
\psi^{-1}\gamma\psi$; hence we obtain $\End(X) \hookrightarrow \End(H)$; this gives a finite index
subgroup $\Aut(X) \to \Aut(H)$; hence the number of $\lambda$ which
given $\zeta = \psi^{\ast}(\lambda)$ is finite.  \B 

\subsection{}\label{defcs}{\bf The central stream.}
For a symmetric Newton polygon $\xi$ we have defined the minimal $p$-divisible group
$H = H(\xi)$. We see that $H$ admits a principal quasi-polarization $\zeta$ 
defined over $\FF_p$, and
we have  seen, see  \ref{unique}, that it is unique up to isomorphism 
over an algebraically closed field.
We write
$$\cZ_{\xi} := \cC_{H(\xi)}(\cA) = \cC_{(H(\xi), \zeta)}(\cA);$$
this will be called the {\bf central stream} in $\cA = \cA_{g,1} \otimes \FF_p$ 
defined by $\xi$.
There does exist $[(X,\zeta)] = x \in \cA$ with $(X,\zeta)[p^{\infty}] \otimes k 
\cong (H(\xi),\zeta)$, such that $\zeta$ is a principal polarization on $A$;
in this case we have $\cZ_{\xi} = C(x)$.

\subsection{}{\bf Remark.} There exist $p$-divisible groups over $k$  for which the number 
of principal quasi-polarizations is bigger than one.\\
{\bf Example.} If $X \sim (G_{1,1})^2$ and $a(X)$ = 1, the number of principal 
quasi-polarizations on 
$X$ over $k$ is bigger than one.

\subsection{}\label{isog} {\bf Proposition} (i)  Type {\bf (I$_r$)}.  
{\it There exists an isogeny $F: (G_{1,1}, \tau_{r+1}) \to (G_{1,1}, \tau_{r})$.}\\
(ii) Type {\bf (II$_r$)}. {\it There exist isogenies 
$\beta: (G_{1,1}, \tau_{r})^2  \to ((G_{1,1})^2,\nu_r)$ and 
$\gamma: ((G_{1,1})^2,\nu_{r+1}) \to (G_{1,1}, \tau_{r})^2$.}\\
(iii) {\it There exists an isogeny 
$(\pi,1): (H_{m,n} \oplus H_{n,m}, \zeta_{r+1}) \to (H_{m,n} \oplus H_{n,m}, \zeta_r)$}.\\
(iv) {\it For quasi-polarized $p$-divisible groups $(X,\lambda)$ and $(Y,\mu)$ with 
$\cN(X) = \xi = \cN(Y)$ over $k$ there exist a quasi-polarization $\zeta$ on $Z:=H(\xi)$
and isogenies}
$$(X,\lambda)  \quad \stackrel{\va}{\longleftarrow} \quad (Z,\zeta)\quad 
\stackrel{\psi}{\longrightarrow}\quad(Y,\mu), 
\quad i.e. \quad  \va^{\ast}(\lambda) = \zeta  = \psi^{\ast}(\mu).$$
{\bf Proof.} Statement (i) follows by a direct verification. \\
(ii) Choosing $\alpha_p \subset (G_{1,1})^2$ with quotient superspecial we construct $\beta$;
direct verification shows the quotient is of type (I$_r$)$^2$; e.g. one can choose 
$\beta = (F,1): (G_{1,1})^2 \to (G_{1,1})^2$. 
Choose $\gamma: X \to (G_{1,1})^2$ of degree $p$ and $X$ superspecial. \\
(iii) This follows by direct verification.\\
(iv) It suffices to show this for the isoclinic parts separately. For the supersingular part 
this follows from \ref{ss} and (i) $\sim$ (iii). The non-supersingular case follows from
(iv) and \ref{eld}.  
\hfill$\Box$

\subsection{}\label{dim=}{\bf Theorem.} (i) {\it  Let $[(A,\lambda,h)] = x \in \cA_{g,d,n} 
\otimes \FF_p$;
here $d \in \ZZ_{>0}$, i.e. we consider polarizations of degree $d^2$, and $n \in \ZZ_{\geq 3}$
is not divisible by $p$; 
suppose $K$ is a perfect field, and suppose that $(X,\lambda) = (A,\lambda)[p^{\infty}]$
is defined over $K$; we view the locally closed set $\cC_{(X, \lambda)}(\cA_{g,d,n} \otimes K)$
as a locally closed subscheme, denoted by $C$,  with induced, reduced scheme structure over $K$. 
For $x, y \in C(k)$ there is an isomorphism of formal schemes
$$C^{/x} \quad \cong \quad  C^{/y}.$$
This scheme $C$
is} smooth over $K$. {\it Any irreducible component of $\cW_{\xi}^0(\cA_{g,d,n} \otimes k)$ which contains $x$ also contains $C$.}\\
(ii) {\bf Notation} $c(-)$. {\it For every symmetric Newton polygon $\xi$ there is 
a number $c(\xi)$ such that 
for every for every $d \in \ZZ_{>0}$, and every $x \in \cC_{(X,\lambda)}(\cA_{g,d} \otimes \FF_p)$  
with $\cN(A) = \xi$, the 
scheme $C$ considered in} (i) {\it 
is pure of dimension $c(\xi)$.}\\
{\bf Proof.} (i) It suffices to prove part (i) in case we work over $K = k$, an 
algebraically closed field. Let $\cD = \Def(X,\lambda)$ be the  formal deformation 
functor. Consider $\cD = \Spf(R)$. We know that the universal $p$-divisible group 
over $\cD$ can be given over $\Spec(R)$, see \cite{dJ}, Lemma 2.4.4 on page 23;
we consider 
$$\cC_{(X,\lambda)}(\Spec(R) = E \subset \Spec(R).$$
Consider $[(A,\lambda)] = x,  [(B,\mu)] = y \in  C(k)$, with  
$$(A,\lambda)[p^{\infty}] = (X,\lambda) =  (B,\mu)[p^{\infty}].$$
By a theorem by Serre and Tate, see \cite{KatzST}, Th. 1.2.1, these identifications
give  isomorphisms
$$\Def(A,\lambda) \quad\stackrel{\sim}{\longleftarrow}\quad  \Def(X,\lambda)  
\quad\stackrel{\sim}{\longrightarrow}\quad  \Def(B,\mu).$$
These induce identifications
$$C^{/x} \quad=\quad  E^{/0}\quad=\quad C^{/y},$$
where $C^{/x}$ denotes the formal completion of $C$ at the closed point $x \in C$.
Hence the reduced scheme $C$ over $k$ is generically regular, hence  
is regular. 

For  $y \in C(k)$ with the indentifications above, we obtain an identification $\cW_{\xi}(\Def(A,\lambda)) \cong \cW_{\xi}(\Def(B,\mu))$. This proves that the number of irreducible components of the formal completion $\cW_{\xi}(\cA_{g,d,n} \otimes k)^{/y}$ for all $y \in C$ is constant. This proves the last statement in (i).

\vn
(ii) It suffices to prove this statement at points over an algebraically closed 
field $k$. By the proof in  (i) it follows that all irreducible components of 
$ \cC_{(X,\lambda)}(\cA_g \otimes k)$ have the same dimension.
Suppose $A[p^{\infty}] = X$ and $(A, \lambda) = x \in \cC_{(X,\lambda)}(\cA_g \otimes k)$, and let 
$(Y,\mu)$  be a quasi-polarized $p$-divisible group with $\cN(X) = \xi = \cN(Y)$.
W are going to show that the dimensions of $\cC_{(X,\lambda)}(\cA_g \otimes k)$, 
respectively $\cC_{(Y,\mu)}(\cA_g \otimes k)$  are the same. (Once this is proven,
we see that this dimension only depends on $\cN(X) = \xi$, and that number will be baptized 
$c(\xi)$.)

By \ref{isog} we construct isogenies of quasi-polarized $p$-divisible groups 
$$(X,\lambda) \stackrel{\va}{\longleftarrow} (Z,\zeta)
\stackrel{\psi}{\longrightarrow}(Y,\mu);$$
we conclude the existence of isogenies of polarized abelian varieties 
$(A,\lambda) \stackrel{\va}{\longleftarrow} (M,\zeta)
\stackrel{\psi}{\longrightarrow}(B,\mu),$ with $(M,\zeta)[p^{\infty}] = (Z,\zeta)$
and $(B,\mu)[p^{\infty}] = (Y,\mu)$. The claim (ii) follows from the following
lemma applied to both $\va$ and $\psi$; the lemma says ``{\it an isogeny between polarized abelian varieties extends to an isogeny correspondence between the two leaves}''.
 
\subsection{}\label{LemT}{\bf Lemma.}  {\it Let 
$(A,\lambda) \stackrel{\va}{\longleftarrow} (M,\zeta)$,
 an isogeny, and  
$$(A,\lambda)[p^{\infty}] = (X,\lambda) \stackrel{\va}{\longleftarrow} (Z,\zeta) = (M,\zeta)[p^{\infty}]$$
 be defined over an algebraically 
closed field $k$; let $\deg(\va) = p^i$. Write $x = [(A,\lambda)] \in \cA_{g,\ast,n} \otimes k$} (in the notation we omit the level structure with $n \geq 3$ prime to $p$), {\it  and let 
$C_x $ be the irreducible component of $\cC_{(X,\lambda)}(\cA)$  containing $x$; write $z = [(M,\zeta)]$  and let $C_z$  be
the irreducible component of $ \cC_{(Z,\zeta)}(\cA)$ containing $z$. There exist an isogeny correspondence consisting of: a scheme
$T$ and finite surjective morphisms $C_x \twoheadleftarrow T \twoheadrightarrow C_z$. 
Hence $\dim(C_x) = \dim(C_z)$.}\\
{\bf Proof.} Consider the universal family restricted to this locally closed subscheme $C_z$; we obtain $(\cM,\zeta) \to C_z$.
Consider the universal family restricted to  $C_x$; we obtain $(\cP,\zeta) \to C_x$. Let $h$ be the height of $Z$; 
choose an integer $N \geq N(h) + i$, where $N(h)$ is as in \ref{2}. By \ref{fconst} we construct a finite, surjective morphism $h: T \twoheadrightarrow C_x \subset C(z)$ with 
an isomorphism $\cM[p^N]_T \cong Z[p^N] \times T$; assume $T$ is irreducible. 
Let $\Ker(\va: \cM_z = Z  \to  X) = G$.
Note that $G$ is isotropic for the form given by $\zeta$.
We obtain $G \times T \subset Z[p^N] \times T \subset \cM_T$.
This gives $(\cX',\lambda') \to T$ by: 
$$\va': \cM_T \to \cM_t/(G \times T) =: \cX', \ \ \  (\va')^{\ast}(\lambda') = \zeta$$
(note that $\zeta$ descends via $\va'$). This family over $T$ plus the induced 
level-$n$-structure from $z$ defines a morphism  $f': T \to \cA_{g,\ast,n} \otimes k$.

Write $[p^{N(h)}]^{-1}(G) =: G' \subset Z[p^N]$. For any $t \in T(k)$ we obtain an isomorphism $G'/G \cong \cX'_t[p^{N(h)}]$. By \ref{2}  we see that  morphism $f'$ factors
through $f: T \to C_x \subset C(x) \subset \cA_{g,\ast,n} \otimes k$. Hence we have {\it constructed}
$$C_x \quad\stackrel{f}{\twoheadleftarrow}\quad T \quad\stackrel{h}{\longrightarrow}\quad C_z.$$ 
{\it We show  that $f$ is  quasi-finite}. It suffices to show that for any point $u \in C_x(k)$ the fiber above $u$ in $f: T \to C_x$ is finite. Suppose $S \subset T$ is an irreducible, reduced subscheme with $f(S) = u$; we show that $S$ is a point.
For any $y \in C_x$ choose an isomorphism $(\cP_x,\lambda)[p^{\infty}] \cong (X,\lambda)$, and for any $t \in C_z$ choose  an isomorphism $(\cM_t,\lambda)[p^{\infty}] \cong (Z,\lambda)$
Choose $q = p^m$ with $m \geq i$, and construct $\Lambda = q^{\ast}(\lambda)$ and
isogenies 
$$(\cP,\lambda)_T \stackrel{\va}{\longleftarrow} (\cM,\zeta)_T \stackrel{\beta}{\longleftarrow} 
(\cP,\Lambda)
\quad\mbox{with}\quad  \va{\cdot}\beta = q.$$
A point $s \in S \subset T$ gives isogenies $\cP_{f(s)} \cong X \to \cM_s \cong Z \to \cX'_s \cong X$. By \ref{5} the number of possible kernels of $\beta_s: \cP_{f(S)} \cong X \to \cM_s \cong Z$ is finite; these kernels are in $\cX_S[q]$; by \ref{4} and by the fact that $S$ is irreducible, the resulting morphism $S \to Gr$, in the notation of \ref{4}, is constant.  Hence $(\cM,\zeta)_S \to S$ is \gfc. As $S \to C_z$ is finite, and as this morphism is the moduli map defined by  $(\cM,\zeta)_S \to S$ we conclude that the image of $S \to C_z$ is finite (hence it is just one point); we conclude that $S$ has dimension zero. We have proved that $f: T \to C_x$ is  quasi-finite. {\it We see that the existence of an isogeny  $\va: (M,\zeta)\longrightarrow(A,\lambda)$ proves} $\dim(C_z) \leq \dim(C_x)$.

Let $x' = [(A,\lambda)]$. Over $C_x$ we have a family $(\cP,\Lambda)$ this defines a morphism $C_x \to \cA_{g,\ast}$ which factors through $C_{x'}$. This map is injective; hence $\dim(C_x) \leq \dim(C_{x'})$. The existence of $\beta: (A,\Lambda) \to (M,\zeta)$ proves, by the arguments above that 
$\dim(C_{x'}) \leq \dim(C_z)$. {\it We conclude that $\dim(C_z) = \dim(C_x)$, and that $f: T \to C_x$ is dominant.}

{\it We are going to prove that $f: T \to C_x$ is 
surjective and finite}. As the image  $f(T)$ is dense in $C_x$, and $f$ is quasi-finite, it suffices to prove the valuative criterion for $f$. We use the notation $q, \beta$ introduced above. Let $K \supset R \to R/m_R=k$ be a discrete valuation ring with field of fractions and residue class field. For any  $y \in C_x$ we choose a commutative diagram

$$\begin{array}{cccc}
\Spec(K) & \stackrel{\eta}{\longrightarrow} & T &\\
\downarrow &  & \downarrow f&\\
\Spec(R) & \stackrel{g}{\longrightarrow} & C_x,& \quad \quad  g(\Spec(k)) = y \in C_x.
\end{array}$$
We obtain $\beta_{\Spec(K)}: (\cP,\Lambda)_{\Spec(K)} \to (\cM,\zeta)_{\Spec(K)}$. We have $\Ker(\beta_{\Spec(K)}) \subset  (\cP,\Lambda)_{\Spec(K)}$; the existence of its flat extension to $(\cP,\Lambda)_{\Spec(R)}$ showns that $(\cP,\Lambda)_{\Spec(K)} \to (\cM,\zeta)_{\Spec(K)}$  extends to  $(\cP,\Lambda)_{\Spec(R)} \to (\cM',\zeta)_{\Spec(R)}$ over $\Spec(R)$. The moduli map defined by $(\cM',\zeta)_{\Spec(R)} \to \Spec(R)$ by \ref{Cp} lands into $C_z$.
 As $T \to C_z$ is finite, hence proper, we conclude that $\eta$ extends to $\Spec(R) \to T$ leaving the diagram commutative. This proves the valuative criterion for $f: T \to C_x$. Hence $g(\Spec(k)) = y \in C_x$ is in the image of $f:T \to  C_x$; we conclude $f$ is surjective. By the valuative criterion, $f$ is proper; hence $f$ is finite. This ends the proof of the lemma. \hfill$\Box$\ref{LemT}

\vn
In the first part of the theorem we have seen that components of $C(x)$ have the same dimension. Using \ref{isog} and \ref{LemT} we see that any two central leaves in the same open Newton polygon stratum admit a finite-to-finite isogeny correspondence. Hence,  in the notation of the theorem, this proves that $\dim(C(x)) = \dim(C(z)) 
= \dim(C(y))$; this concludes the proof of the theorem.
\hfill$\Box$\ref{dim=}

\vn
Smoothness as in the theorem of a central leaf was already known in a special case
in a different setting, see 
\cite{HT}, Coroll. III.4.4 on page 114.

\subsection{}
We see the curious phenomenon:  

{\it although the dimension of a component of the Newton 
polygon strata depends on the degree of the polarization} (and there are Newton polygons
for which this dimension is not constant inside that stratum, compare \ref{maxdim}), 

{\it the dimension of the central leaves in
one Newton polygon stratum in $\cA_g \otimes \FF_p$ are all of the same dimension.} 

\vn
Moreover, {\it in general a Newton polygon stratum is singular} (even after adding level 
structure, even after considering only  one irreducible component of a Newton polygon 
stratum); however

{\it every central leaf} (level structure considered) {\it is smooth over the base field.}

\vn
Isogeny correspondences in general blow up and down (in characteristic $p$ if local-local 
kernels are involved); however

{\it isogeny correspondences between central leaves are finite-to-finite.}

\subsection{}\label{icorr}{\bf Remark.} Using \ref{LemT} and methods as in \cite{FC}, VII.3 we can show:\\
{\it Suppose $\va: (A,\lambda) \cdots\to (B,\mu)$ is a quasi-isogeny, then there exists an isogeny correspondence
$$C_x \quad\stackrel{f}{\twoheadleftarrow}\quad T \quad\stackrel{h}{\twoheadrightarrow}\quad C_z$$
with $x = [(A,\lambda)]$ and $ z = [(B,\mu)]$ which contains the isogeny correspondence $\va$ and such that $f$ and $h$ are finite and surjective.}

\subsection{}\label{cxi}{\bf Remark.} In \cite{C.FO} we prove: 
for a symmetric Newton polygon $\xi$ of height $h = 2g$ we have:
$$c(\xi) = 2{\cdot}\sum_{0 < j \leq g} \ \  (\sigma(j) - \xi(j)).$$

\section{Isogeny leaves}
\n
Isogeny leaves, as defined below, can be considered as in \cite{RZ} in the context of ``Rapoport-Zink spaces''. We have chosen to give an independent exposition.

\subsection{}\label{defH}{\bf Definition.} {\it Let $S$ be a scheme over $k$, 
and let $(\cD, \lambda) \to S$ be a polarized abelian scheme over $S$. 
A  reduced subscheme $I \subset S$ 
is called an $H_{\alpha}$-subscheme in $S$ if there exists a polarized 
abelian variety $(M,\zeta)$
over $k$, a scheme $T$ of finite type over $k$, a surjective morphism 
$T \twoheadrightarrow I$ and an isogeny 
$\va: (M,\zeta) \times T \to (\cD, \lambda) \times_I T$ such that 
every geometric fiber of $\va$ is a successive extension of subfactors 
isomorphic with $\alpha_p$;} equivalently: every fiber of $\Ker(\va)$ 
is of local-local type. 
{\it The pair $(T \to I, \va)$ will be called a chart for the $H_{\alpha}$-scheme.}

\subsection{}\label{isogl}{\bf Theorem}
($I(x)$: the existence of certain maximal $H_{\alpha}$-subschemes).   
{\it Consider $d \in \ZZ_{\geq 1}$; we write $\cL:= \cA_{g,d} \otimes k$. 
Let $x \in \cL(k)$. 
There exists a reduced, closed subscheme
$I(x) \subset \cL$ which is $H_{\alpha}$ such that every irreducible component of 
$I(x)$ contains $x$, and such that for every integral subscheme 
$V \subset \cL$ we have:
$$V \subset \cL  \quad\mbox{is an } \ H_{\alpha}-\mbox{subscheme, and} \quad x \in V 
\quad\Longrightarrow\quad   V \subset I(x),$$}
i.e. $I(x)$ is the union of all irreducible $H_{\alpha}$-sets in $\cL$ containing $x$.

\subsection{}\label{defil}{\bf Notation.} Let $[(A,\lambda)]=x \in \cA_g(k)$; let 
$I(x) \subset \cA_{g,d} \otimes k$ be as in the 
theorem. An irreducible component of $I(x)$ with the induced, reduced scheme structure 
will be called an {\it isogeny leaf}. I.e. {\it an isogeny leaf is a maximal,  
integral $H_{\alpha}$-subscheme of $\cA_g \otimes k$.}

\subsection{}{\bf Corollary.} {\it Let $S$ be a scheme over $k$, and $(\cD, \lambda) \to S$ 
be a polarized abelian scheme; let $x \in S(k)$. 
There exist a closed subscheme $I(x) \subset S$ which 
is $H_{\alpha}$ and which contains all irreducible $H_{\alpha}$-subschemes $V$ such that 
$x \in V \subset S$.}    \B

\subsection{}\label{IsogCorr} {\bf Isogeny correspondences.} Consider the scheme of all isogeny correspondences 
$$\cI \longrightarrow \left(\cA_{g,\ast} \otimes \FF_p\right) \times \left(\cA_{g,\ast} 
\otimes \FF_p\right).$$ This scheme is a coarse moduli scheme for all diagrams of isogenies
$$(A,\lambda) \stackrel{\va}{\longleftarrow} (M,\zeta)
\stackrel{\psi}{\longrightarrow}(B,\mu).$$
See \cite{FC}, page 251 and the rest of VII.3 and VII.4 in \cite{FC} for a discussion.
We shall also use this correspondence in case of isogenies between polarized abelian
varieties with a level-$n$-structure, with $n$ not divisible by $p$, and considering 
isogenies of degree prime to $n$.
Note:

\subsection{}\label{Closed}
{\bf Lemma.} {\it Let $G \to S$ be a finite, locally free group scheme of rank $p^n$. The set of points 
$s \in S$ such that 
$G_s$ is a local group scheme is closed.}\\
This is a question in characteristic $p$. The condition $G_s \subset G[F^n]$ inside the flat scheme
$G \to S$ is a closed condition. \B 

\vn
We define $\cH_{\alpha} \subset \cI$ to be the  reduced subscheme defined 
by the condition 
that we consider only isogenies, in the notation above: $\va$ and $\psi$, 
with local-local kernel. By the previous lemma we see that $\cH_{\alpha} \subset \cI$
is a closed subscheme.

\subsection{}{\bf Corollary.} {\it Let $(T \twoheadrightarrow I, \va)$ be a chart as 
in the definition} \ref{defH};  
{\it there exists a factorization  $T \to T' \to I$, and a morphism of charts 
$(T, \va) \to (T', \va')$ for $I$ such that $T' \twoheadrightarrow I$ 
is a proper morphism (of finite type). One can choose $T' \twoheadrightarrow I$ 
to be finite above all generic points of $I$.

If $I^0 \subset S$ is $H_{\alpha}$ as in the definition, then its closure 
$I \subset S$ is also $H_{\alpha}$.}\\
Every component of $\cH_{\alpha}$ is of finite type and proper over both factors.  \B

\subsection{}{\bf Remark.} 
In general, a connected union of  $H_{\alpha}$-subschemes
need not be a closed subscheme.\\
{\bf Example.} Let $\xi = (2,1) + (1,2)$, and consider the central stream $\cZ_{\xi}$.
Let $x \in \cZ_{\xi}$. Consider all $x'$ in $\cA = \cA_{3,1} \otimes \FF_p$ 
obtained over some field by iterated $\alpha_p$-isogenies from $x$. 
One shows that this gives a connected union of $H_{\alpha}$-subschemes; this union is
finite, i.e. it is a closed subset,  iff $x \in \cA(\overline{\FF_p})$. 
So, in order to produce a non-closed, connected 
union of $H_{\alpha}$-subschemes, start with $x \in \cZ_{\xi}$ which is not 
defined over a finite field.

\subsection{}\label{constant}{\bf Lemma.} {\it Let $T$ be a connected, integral $k$-scheme; 
let $Y$ be a 
$p$-divisible group over $k$,  and let $\cZ \to T$ be a $p$-divisible group which is \gfc. 
Assume there exists an isogeny $\psi: Y \times T \to \cZ$. Then $\cZ \to T$ is constant.}\\
{\bf Proof.} The finite flat group scheme $\cG = \Ker(\psi) \to T$, say of rank $p^i$, 
defines a section in $Gr_{Y,i} \times T \to T$, see \ref{4}. By \ref{5} this section lands 
in every fiber in a finite set  of closed points $\Psi \subset Gr_{Y,i}$, see \ref{5}. As $T$ is connected
this section is constant. Hence there exists $G \subset Y$ such that $\cG = G \times T$.
This shows that $\cZ \cong (Y/G) \times T$.  \B    

\subsection{}\label{Pfisogl}{\bf Proof of} \ref{isogl}. We choose $n \in  \ZZ_{\geq 3}$
not divisible by $p$. Let $\xi$ be the Newton polygon of $x$.\\
{\bf Claim.} {\it If $V \subset \cW_{\xi}^0(\cA_{g,d,n} \otimes k)$ is closed and $H_{\alpha}$, then $V$ is closed and $H_{\alpha}$ in  $\cA_{g,d,n} \otimes k$ .}\\
Let $T' \supset V$ be the Zariski closure of $V$ in  $\cA_{g,d,n} \otimes k$,  let $T' \to T$ be the normalization, and let $V'\subset T'$   be the inverse image of $V$. We find a polarized abelian variety $(M,\zeta)$ and an isogeny $\va': (M,\zeta) \times V' \to (\cD,\lambda) \times_V V'$ such that every geometric fiber of $\Ker(\va') = G'$ is of local-local type. By \cite{FC}, Proposition I.2.7 there exists an extension $\va: (M,\zeta) \times T' \to (\cD,\lambda) \times_T T'$; it follows that $\Ker(\va)$ is of local-local type. Hence the image of $T$ in $\cA_{g,d,n}$ is contained in $\cW_{\xi}^0(\cA_{g,d,n})  \otimes k$; i.e. $V$ is closed and $H_{\alpha}$ in $\cA_{g,d,n} \otimes k$, which proves the claim.

\vn
Write $x = (A, \lambda, f)$ and let $\xi = \cN(A)$. Let
$W_1, \cdots , W_s$ be the irreducible components of $\cW_{\xi}^0(\cA_{g,d,n} \otimes k)$ 
which contain the point $x$.

\vn 
{\bf Maximal $H_{\alpha}$-schemes.} {\it For every $W_j$ we construct
a closed subscheme   $I'_j(x) \subset W_j \subset\cA_{g,d,n} \otimes k$ which is $H_{\alpha}$, 
such that every irreducible component of $I'_j(x)$ contains $x$  and such that
every $x \in V \subset W_j$ which is irreducible and $H_{\alpha}$ is contained in $I'_j(x)$.}

\vn
We choose some $1 \leq j \leq t$ and we write $W' = W_j$.
Let $W \to W'$ be the normalization map. Choose $y \in  W(k)$
with $y \mapsto x \in W'$.
Let $(\cU', \lambda, f) \to \cL$ be the universal family; write  
$(\cU, \lambda, f) = (\cU', \lambda, f) \times _{\cL} W$. 

\vn
{\it We show that there exist a polarized abelian scheme $(\cM, \zeta,f)$ with level-n-structure
over $W$ and an isogeny $\va: (\cM, \zeta,f) \to (\cU, \lambda, f)$ with local-local kernel over $W$ 
such that $\cZ := \cM[p^{\infty}]$  is \csd.}  -- In fact this follows from
\cite{FO.Z}, 2.1 applied to $\cU[p^{\infty}]$: as  $W$ is normal there exists an isogeny 
$\va: \cZ \to \cX :=  \cU[p^{\infty}]$ where $\cZ \to W$ is a \csd \ $p$-divisible group;  
this isogeny 
defines $\va: (\cM, \zeta, f) \to (\cU, \lambda, f)$ as desired. 

\vn
Let $\deg(\va) = p^i = q$ (remark: this degree
may depend on the choice of the component $W_j$). We write $Z = \cZ_y$ and 
$(M,\zeta) = (\cM, \zeta)_y$; hence $M[p^{\infty}] = Z$. 
We write $[(M,\zeta,f)] =:m \in \cA_{g,dq,n}(k)$.  Note that $(\cM, \zeta, f) \to W$
defines a correspondence 
$$\cA_{g,dq,n}\otimes k \quad\leftarrow\quad \cJ \quad\rightarrow\quad W; 
\quad\mbox{define}\quad I_j(x) = \cJ^{-1}(m) \subset W;$$
here $\cJ^{-1}(m)$ is the projction to $W$ of the inverse image
of $m$ in $\cJ$.
Define $x \in I'_j(x) \subset M' = W_j$ to be its image in $W'$.

\vn
Claim: {\it For any irreducible $V \subset W_j = W' \subset \cL$, which is $H_{\alpha}$,
with $x \in V$  
we have $V \subset I'_j(x)$.} -- Indeed,  
consider a chart $(T' \twoheadrightarrow V, \rho)$ of finite type
for $V$, with $\rho: D \times T' \to \cU'_{T'}$. Let $T$ be an irreducible component
of $T' \times_{W'} W$ mapping onto $V$ which contains a point  $t \in T$ with $t \mapsto y \in W$. 
We have an isogeny 
$\rho_T: D \times T \to \cU'_T = \cU_T$, and we have an isogeny 
$(\cM \to \cU) \times_W T = \cM_T \to \cU_T$. By \ref{constant} we conclude that $\cM_T \to T$
is constant; in fact, $\cM_T \cong M \times T$; the polarizations $\zeta_T$ and 
$\zeta \times T$ coincide at $t$, hence they are equal; thus we have an isogeny
$\va_T: (\cM_T,\zeta_T) \to (\cU_T,\lambda_T)$. The kernel of the isogeny 
$\va_T: \cM_T \cong M \times T \to \cU_T$ we denote by $\cG' \to T$; this defines a morphism
$f: T \to Q$ such that $f^{\ast}(\cG) = \cG'$. As this is the kernel of an isogeny of 
polarized abelian schemes we conclude that $f$ factors: 
$f: T \to P \subset Q$. The commutative diagram
$$\begin{array}{ccc}
T & \rightarrow & P\\
\downarrow &  & \downarrow\\
V & \rightarrow & W'
\end{array}$$
shows that
  the image is contained in $I_j(x)$: 
$$V = {\rm Im}\left(T \to P \stackrel{\Psi}{\longrightarrow} W'\right) \subset I_j(x).$$
Hence $I'_j(x)$ thus constructed contains all irreducible $H_{\alpha}$-subschemes
contained in $W' = W_j$ passing through $x$.  
Hence $I(x) := \cup_j \  I'_j(x) \subset \cL \subset \cA_g \otimes k$
is $H_{\alpha}$, closed in  $\cA_g \otimes k$,
and it contains all irreducible  $H_{\alpha}$-subschemes in  $\cA_g \otimes k$
containing $x$. This proves the theorem.         

\B \ref{isogl}

\subsection{}{\bf  Proposition.}\label{compl}  {\it An irreducible, maximal  $H_{\alpha}$-scheme  $I \subset \cA_g \otimes k$ is closed in  $\cA_g \otimes k$ and proper over k.}
\\
{\bf Proof.} This follows from  the definition of maximal  $H_{\alpha}$-scheme, and from \ref{4}, using the valuative criterion of properness.  \B

\subsection{}{\bf Proposition.}\label{dim=i} {\it Let $x \in \cA_g(k)$ with $x = [(A,\lambda)]$,
and $(X, \lambda) := (A, \lambda)[p^{\infty}]$. 
Let  $[(Y, \mu)] =:y \in C(x)(k) = \cC_{(X, \lambda)}(\cA_g \otimes k)(k)$.
Then there is an isomorphism of formal schemes} $$I(x)^{/x} \cong I(y)^{/y}.$$\\
{\bf Proof.} As $y \in C(x)(k)$ there exists an isomorphism
$$\va: (X, \lambda) \quad\stackrel{\sim}{\longrightarrow}\quad (Y, \mu).$$
This gives $\va: D_X := \Def(X, \lambda) \stackrel{\sim}{\longrightarrow} D_Y :=\Def(Y, \mu)$, and by the Serre - Tate theorem we obtain $\va: D_A :=\Def(A, \lambda) \stackrel{\sim}{\longrightarrow} D_B:=\Def(B, \mu)$. \\
{\bf Claim.} Using these isomorphisms, 
{\it from  the facts that 
$I(x)$ is $H_{\alpha}$ and  $I(y)$ is maximally $H_{\alpha}$ at $y$ we can conclude that $\va(I(x)^{/x}) \subset I(y)^{/y}$ under the given identifications.}

After considering an appropriate level structure, we conclude we have an isogeny $(M,\zeta) \times I(x) \to (\cB,\lambda)$ expressing that $I(x)$ is $H_{\alpha}$. Let the degree of this isogeny be $p^i$. Considering the $p$-divisible groups involved, and restricting to $I(x)^{/x}$ we obtain an isogeny $(M,\zeta) \times I(x)^{/x} \to (\cD,\lambda) \times_{D_{A}} I(x)^{/x}$. Consider $\va(I(x)^{/x}) \subset I(y)^{/y} \subset D_Y$. Using $\va: (X, \lambda) \stackrel{\sim}{\longrightarrow} (Y, \mu)$ we obtain an isogeny  $(M[p^{\infty}],\zeta) \times \va(I(x)^{/x}) \to (\cY,\mu) \times_{D_Y} \va(I(x)^{/x})$.

Consider the isogeny correspondence, see \ref{IsogCorr}, restricted to the factors considered: 
$$\cI \subset \left(\cA_{g,d,n} \otimes k\right) \times  \left(\cA_{g,id,n} \otimes k\right).$$
Taking the projection on the second factor to be equal to  $(M,\zeta)$ (we omit the level structure in the notation), and only local-local kernels in the isogeny, we obtain a closed set in $\cI$, see \ref{Closed}, and using the first projection (which is a proper morphism) we obtain a closed set $T \subset \cA_{g,d,n} \otimes k$. By construction it contains the point $y$, it is $H_{\alpha}$, and the completion of $T$ at $y$ contains $\va(I(x)^{/x})$. Because $I(y)$ is maximally $H_{\alpha}$ at $y$ we see that every reduced, irreducible component of $T$ containing $y$ is contained in $I(y)$:
$$T_{\rm red} \subset I(y).$$ 
Hence  
$$\va(I(x)^{/x}) \subset (T_{\rm red})^{/y} \subset I(y)^{/y}.$$
This proves the claim. Reversing the roles of $x$ and $y$ we show that $\va^{-1}(I(y)^{/y}) \subset  I(x)^{/x}$. this finishes the proof of the proposition.
  \B

\subsection{}{\bf Remark.} We did not give a relevant definition of $H_{\alpha}$-schemes for families $p$-divisible groups; we did not define 
isogeny leaves in such families.
I don't know in which generality  an equivalent of \ref{isogl} for $p$-divisible 
groups  holds.

\section{The product structure defined by central and isogeny leaves}
\subsection{}\label{zero} {\bf Proposition.}  {\it Let $x \in \cA_g(k)$. 
Every irreducible component of $C(x) \cap I(x)$ has dimension equal to zero.}\\
{\bf Proof.} We will assume considering everything inside $\cL := \cA_{g,d,n} \otimes k$, 
with $n \in \ZZ_{\geq 3}$ prime to $p$. Let $(\cU,\mu,f) \to \cL$ be the universal family.
Let $S$ be a reduced, irreducible component 
of $C(x) \cap I(x)$. 
There exists a chart of finite type 
$(T \twoheadrightarrow S, \va: (M,\zeta) \times T \to (\cU,\mu) \times_{\cL} T))$.
Suppose $T$ is ireducible.
As $S \subset C(x)$, we conclude that $\cU_S$ is \gfc; hence  $(\cU,\mu) \times_{\cL} T$
is \gfc.  Using \ref{4} and \ref{5} we conclude that $(\cU,\mu,f) \times_{\cL} T$ is constant.
Hence the induced moduli morphism $T \twoheadrightarrow S \subset \cL$ is constant. This proves
the proposition.   \B

\subsection{}\label{P}{\bf Lemma.} {\it Let $(Z,\zeta)$ be a quasi-polarized 
$p$-divisible group 
over a field $K$, and let $i \in \ZZ_{>0}$; write $\cK = \Ker(\zeta)$. 
Consider $Gr = Gr_{Z,i}$ as in} \ref{4} {\it with the 
universal family $\cG \to Gr$. The set
$$ P := \{t \in  Gr \mid \cG_t \subset \Ker(\zeta), \mbox{  and it is isotropic for this form}\}$$
is  a closed subset of   $Gr$. }\\
{\bf Poof.} The property 
$\cG_t \subset \cK \times S \subset Z \times S$ is a closed condition 
(flat extensions inside a flat  $Z \times S \to S$);  also the condition that 
$\cG_S^{\bot} \subset \cK/\cG_S \subset Z^t \times S$ is closed; 
this proves the claim.  \B

\vn
We show that components of any open Newton polygon stratum (almost) have a 
product structure 
given by central laves and isogeny leaves contained in that stratum:

\subsection{}\label{TI}  {\bf Theorem} (``central and isogeny leaves  almost give a 
product structure 
on an irreducible component of a
Newton polygon stratum''). {\it Let $d \in \ZZ_{\geq 1}$, let $\xi$ be a symmetric Newton polygon, and let 
$W'' \subset \cA_{g,d} \otimes k$
be an irreducible component of the open Newton polygon stratum 
$\cW_{\xi}^0\left(\cA_{g,d} \otimes k\right)$. There exist integral schemes $T$
and $J$ of finite type over $k$, and a finite surjective $k$-morphism
$$\Phi: \quad T \times J \quad\twoheadrightarrow\quad W'' \subset \cA_{g,d} \otimes k$$
such that 
$$\forall  u \in J(k), \quad \Phi(T \times \{u\}) \quad\mbox{is a central leaf in} \quad W'',$$
every central leaf in $W''$ can be obtained in this way, 
$$\forall t \in T(k),  \quad \Phi(\{t\} \times J)  
\quad\mbox{is  an isogeny leaf in } \quad W'',$$
and every isogeny leaf in $W''$ can be obtained in this way.}\\
{\bf Proof.} Let $W''$ be an irreducible component of $\cW_{\xi}(\cA_g \otimes k)$; choose $n \geq 3$ prime to $p$, let  $$W'' \quad\leftarrow\quad W' \subset \cW_{\xi}(\cA_{g,d,n} \otimes k) \subset \cA_{g,d,n} \otimes k$$ 
be an irreducible (reduced) component of this Newton polygon stratum mapping onto $W''$. The finite morphism $\rho: W \to W'$ is the normalization map. Let $(\cU,\lambda,f) \to W$ be the pull back of the universal family. By \cite{FO.Z}, 2.2,  there exists $\va: (\cM,\zeta,f) \to (\cU,\lambda,f)$ over $W$
such that $\cZ :=\cM[p^{\infty}] \to W$ \csd; write $q = p^i = \deg(\va)$ and $\cG = \Ker(\va)$.
From $(\cM,\zeta,f) \to W$ we obtain the moduli morphism 
$\pi': W \to \cA_{g,dq,n} \otimes k$.

\subsection{}\label{C1}{\bf Claim / Notation.} (i)  {\it The image of the moduli morphism $\pi': W \to \cA_{g,dq,n} \otimes k$ is contained in the central leaf $C  \subset  \cA_{g,dq,n} \otimes k$ passing through any point of  $\pi'(W)$.   We obtain:} 
$$ \cA_{g,d,n} \otimes k \supset W' \stackrel{\rho}{\longleftarrow} \quad W 
\quad\stackrel{\pi}{\longrightarrow}\quad C  
\subset  \cA_{g,dq,n} \otimes k.$$
(ii) {\it The morphism $\pi: W \to C$ is proper.}\\
(iii) {\it For any $w \in W(k)$ and any   irreducible (reduced) component $J$ of the reduced fiber $\left(W \times_C \{\pi(w)\}\right)_{\rm red} \subset W$
the image $\rho(J) \subset W'$ is an isogeny leaf.}\\
(iv) {\it Let $C_{x'} \subset W'$ be a central leaf, and $D$ an irreducible component of $\rho^{-1}(C_{x'}) \subset W$. The morphisms 
$$C_{x'} \twoheadleftarrow D \twoheadrightarrow C$$ 
induced by $\rho$ and by 
$\pi$ are surjective and finite.}\\
(v) {\it The morphism $\pi: W \twoheadrightarrow C$ is surjective.}

\vn
As $W$ is irreducible, and  $\cZ :=\cM[p^{\infty}] \to W$ is \csd, using \ref{3}, we conclude that $\pi'(W)$ is contained in a central leaf. Note that different central leaves have an empty intersection. Hence there is precisely one central leaf $C$ containing  $\pi'(W)$; the restriction of $\pi'$  we call $\pi: W \to C$. This proves (i). 

\vn
Let $R$ be a discrete valuation ring with field of fractions $L$; write $\Gamma = \Spec(R)$ and  $\Gamma^0 = \Spec(L)$;  consider a morphism $f: \Gamma \to C$ such that $f(\Gamma^0) \subset \pi(W)$. We show that also $f(\Gamma) \subset \pi(W)$. Indeed, there exists a finite morphism $g: \Delta \to \Gamma$, write $\Delta^0 = g^{-1}(\Gamma^0)$,  and a morphism $h^0: \Delta^0 \to W$ with $\pi{\cdot}h^0 = (\Delta^0 \to \Gamma \to C)$.  Denote by $(\cM', \zeta') \to \Delta$ the pull back of the universal family over $\cA_{g,dq,n} \otimes k.$ The restriction of this to  $\Delta^0$ is the pull back of $(\cM,\zeta) \to W$.   By \ref{fconst} we can suppose that $\Delta$ is choosen in such a way that $\cM'[p^i]$ is constant over $\Delta$. By \ref{5} and \ref{2} we see that we obtain $(\cM',\zeta) \to (\cU',\lambda)$ over $\Delta$, that the moduli map $\Delta^0 \to W'$ extends, hence $\Delta^0 \to W$ extends to $\Delta \to W$, hence $f(\Gamma) \subset \pi(W)$. Using the valuative criterion, see \cite{HAG} Th. II.4.7 and Exc. II.4.11, this shows (ii), that  $\pi: W \to C$ is proper.

\vn
Clearly $I := \rho(J)$ is $H_{\alpha}$ and irreducible. Let $I'$ be an isogeny leaf containing $I$; suppose $I \subsetneqq I'$; then the $\dim(I') > \dim(I)$. In that case there exists $J \subset J'$ with $\rho(J') = I'$ and $J'$ irreducible; we would obtain $\dim(\pi(J')) > 0$, moreover $\pi(J')$ is $H_{\alpha}$ and $\pi(J') \subset C$. This is a contradition with \ref{zero}. This proves $I =I'$, which proves (iii).

\vn
As $W' \twoheadleftarrow W$ is finite surjective the same holds for $C_x \twoheadleftarrow D$. Note that $D \to C$ is quasi-finite: any positive dimensional $D' \subset D$ with $\dim(\pi(D'))=0$ would give a positive dimensional $\rho(D') \subset C_{x'}$ which is $H_{\alpha}$, a contradiction with \ref{zero}. By (ii) we see that  $D \to C$ is proper, and hence finite;  $\pi(D) \subset C$ is closed. Moreover $\dim(C_{x'}) = \dim(C)$ by \ref{LemT}. Hence $\dim(\pi(D)) =  \dim(C)$. This proves that $\pi: D \twoheadrightarrow C$ is surjective; this proves (iv). As $\pi(D) \subset \pi(W) \subset C$ we have proved (v), and the proof of the claim is finished.          
  \B \ref{C1}
 
\subsection{}{ \bf Construction of $\Phi$.}
We write $\rho'' = \left(W \stackrel{\rho}{\longrightarrow} W' \to W''\right)$.
Choose a point $w \in W(k)$ such that $x'':=\rho''(w) \in W''$   is not contained in any other component of $\cW_{\xi}^0\left(\cA_{g,d} \otimes k\right)$; write $\rho(w) = x'$ and $\pi(w) = z = \pi(J)$; choose $J \subset W$ as in  (iii) of \ref{C1}  with moreover $w \in J$.  We write $M = \cM_w = \cM'_z$, and $Z = M[p^{\infty}]$.  We see that there exists a finite group
group scheme, $\cG \hookrightarrow M[p^i] \times J$
isotropic under the form given by the polarization $\zeta$ 
 such that $\zeta$ descends to a polarization, say $\mu'$, on  
$(M \times J)/\cG$, and such that 
$$((M \times J)/\cG,\mu') \cong (\cU,\mu)_J.$$ 
The
finite group scheme $\cG \hookrightarrow M[p^i] \times J$ defines a morphism
$e: J \to P \subset Gr$, in the notation of \ref{P}; {\it we claim that the morphism $e$ is quasi-finite}; indeed, let $V \subset J$ be a   closed irreducible subset such that $e(V)$ is a point; this implies that the restriction of $\cG$ to $V$ is constant; hence $(\cU,\mu)_V$ is constant; the fact that the family $(\cU,\mu)$ is universal, and the fact that $W \to W'$ is finite implies that $V$ is finite; this proves the claim. 

\vn
Choose $D \subset W$ with $\rho(D) = C_{x'}$ and $w \in D$ as in (iv) of \ref{C1}.
Using \ref{fconst} we choose a commutative diagram
$$
\begin{array}{rcl}
T&\stackrel{\psi}{\longrightarrow}& S\\
g\downarrow &&\downarrow h\\
C_{x'} \leftarrow  D  &\stackrel{\pi}{\longrightarrow}& C = C_z
\end{array}
$$
such that $g^{\ast}(\cM)[p^i] \cong M[p^i] \times T $, an isomorphism  compatible
with  $h^{\ast}(\cM')[p^i] \cong M[p^i] \times S $, and such that $T$ is normal
and integral. Using the notation $(\cM',\zeta',f')$ for the universal famiily over $\cA_{g,dq.n} \otimes k$ we have: $g^{\ast}(\cM,\zeta,f) = (h{\cdot}\psi)^{\ast}(\cM',\zeta',f')$.

\vn
{}From 
$$\cG \hookrightarrow M[p^i] \times J,  \quad\mbox{and}\quad  M[p^i] 
\times T \hookrightarrow \cM_T, $$ 
the pull back of the universal family over 
$Gr$ to $J \to P \subset Gr$ with $P$ as in \ref{P}, respectively
 the pull back of $\cM \to W$ 
to $T$ via $T \to D \subset W$,  we obtain
$$\cG \times T \quad\hookrightarrow\quad M[p^i] \times T \times J 
\quad\hookrightarrow\quad \cM_T \times J  \quad\mbox{over}\quad T \times J.$$ 
By construction this subgroup scheme is  contained in the kernel of and 
isotropic under the form given by $\zeta$; hence
we can descend the polarization to the quotient, obtaining
$$(\cM_T \times J,\zeta_T \times J,f) \longrightarrow (\cQ,\mu,f),  
\quad \cQ:= (\cM_T \times J)/(\cG \times T) $$
over $T \times J$. 
This defines the moduli morphism:
$$ [(\cQ,\mu,f)] = \Psi': \quad T \times J \quad\longrightarrow\quad \cA_{g,d,n} \otimes k.$$
Write $\Phi = (W'' \leftarrow W'){\cdot}\Psi'$.

\subsection{}\label{C2} 
{\bf Claim / Notation.}  {\it  
\begin{description}
\item[(i)]  We conclude that $\Psi': T \times J \to W'$; 
this defines a unique $\Psi: T \times J \to W$ such that 
$\rho{\cdot}\Psi = \Psi'$; we obtain a commutative diagram: 
$$\begin{array}{ccccccl}
\ \ T \times J & = &  T \times J & = & \ \ T \times J &\stackrel{\rm proj_1}{\longrightarrow}&  T \\
\Phi \downarrow && \Psi' \downarrow && \Psi \downarrow &&\downarrow h{\cdot}\psi\\
 \ \ \ W'' &  \longleftarrow &  \ \ \  W'&   \stackrel{\rho}{\longleftarrow}&
 \ \ \ W &\stackrel{\pi}{\longrightarrow}& C;
\end{array}
$$
moreover, for every $u \in J(k)$ we have:
$$\left(T \cong T\times \{u\} \hookrightarrow T \times J  \stackrel{\Psi}{\longrightarrow} W
\right) = \left( T \stackrel{g}{\longrightarrow} D \hookrightarrow W \right);$$
\item[(ii)] the morphism $\Psi: T \times J \twoheadrightarrow W$ is finite and surjective; hence $\Phi$ is finite and surjective;
\item[(iii)]
 if $u \in J(k)$ then  $\Phi(T \times \{u\})$
is the central leaf passing through $\rho''(\Psi(t,u)) = \Phi(t,u)$ for any $t \in T(k)$;
every central leaf can be obtained in this way;
\item[(iv)]  if $t \in T(k)$ then  $\Phi(\{t\} \times J)$ is an isogeny leaf in $W''$; every isogeny leaf can be obtained in this way. 
\end{description}
}
\n
Choose $t_0 \in T$ with $g(t_0) = w$; choose $u_0 \in J$ with 
$e(u_0) = [\cG_w] \in  P \subset Gr$. By construction we see that
$\Psi'(t_0,u_0) = x' \subset W'$. As the Newton polygon is constant on 
$\Phi(T \times J)$,
and by the choice of $w$ we conclude that  $\Psi'(T \times J) \subset W'$. As $T$ and $J$ 
are normal, and $\rho: W \to W'$ is the normalization morphism, we obtain $\Psi$ 
as desired. As $\rho$ is biregular at $w \in W \to W'$ we see  that $\Psi(t_0,u_0) = w$. Consider the moduli map 
$$\pi' = [(\cM_T \times J,\zeta_T \times J,f)]: T \times J \to C;$$ for every $t \in T(k)$ this morphism restricted to $\{t\} \times J$ gives one point in $C$; hence it factors over ${\rm proj}_1: T \times J \to T$; for every $u \in J(k)$ the restriction of $(\cM_T \times J,\zeta_T \times J,f)$ to $T \times \{u\}$ equals $(\cM_T ,\zeta_T ,f)$, hence $\pi'$  restricted to $T \times \{u\} \cong T$ equals $h{\cdot}\psi$. Hence
$$\pi' = h{\cdot}\psi{\cdot}{\rm proj}_1.$$
We have proved: $(\pi{\cdot}\Psi)^{\ast}(\cM',\zeta') = \Psi^{\ast}(\cM,\zeta) = (\cM_T \times J, \zeta_T \times J)$. This proves that the restriction of $\Psi'$ to $T \times \{w\}$ is given by as the moduli map by this family;  note that $W' \leftarrow W$ is biregular at $w$  by the choice of $w$; hence the last statement in (i) is proved. 
This shows (i).

\vn
{\it We claim that $\Psi$ is quasi-finite}; in fact, let $V' \subset T \times J$  be an irreducible,  closed subset, such that  $\Psi(V')$ is one point; then it is mapped to one point under ${\rm proj}_1 = t \in T(k)$; hence it is  closed in $\{t\} \times J$, i.e. $V' = \{t\} \times V$; this  means that $e(V)$ is one point, where $e: J \to P$ is the morphism considered above ; we have seen that this implies that $V$ is finite;  we conclude that $\Psi$ is quasi-finite. 

By dimension arguments we conclude that 
$\Psi$ is dominant (and we see that $\dim(J) = \dim(W) - \dim(C)$). As $T \to C$ is finite (hence proper), and $J$ is complete, see \ref{compl},  the composition
$$\left(T \times J \to W \to C\right) \quad=\quad \left(T \times J \to T \to C\right)$$ 
is proper.  We conclude by \cite{HAG}, II.4.8 (e) that $T \times J \to W$ is proper. Hence $\Psi: T \times J \to W$ is finite and surjective; this shows (ii).

\vn
Consider $u \in J(k)$; the restriction of $(\cU,\mu,f)$ to  $D_u :=\Psi(T \times \{u\})$
is \gfc; hence  $\rho(D_u)$ is contained in the central leaf $C' := C_{\rho(\Psi(t_0,u))}$ for any $t_0 \in T(k)$. This image is closed in $C'$, and for dimension reasons, see \ref{dim=} , $\rho(D_u)$  is dense  in $C'$; hence equality;  hence $\Psi'(T \times \{u\})$ is a central leaf.

For any $y' \in W'$, choose an irreducible component $D' \subset \rho^{-1}(C_{y'})$, and $w' \in D'$; as  $\Psi: T \times J \to W$ is surjective, we can choose $(t_0,u_0) \in T \times J$ with $\Psi(t_0,u_0) = w'$. We see that $\Psi(T \times \{u_0\}) = D'$; hence the central leaf $C_{y'}$ equals $\Psi'(T \times \{u_0\})$. This proves (iii).

\vn
We see that $\pi(\Psi(\{t\} \times J)) = z'$ is a point. The existence of an irreducible component of $W \times_C \{z'\}$ containing $\Psi(\{t\} \times J)$ of strictly bigger dimension would give a contradiction with (ii) of \ref{C2}, the fact that $T \to C$ is finite, and \ref{zero}. By (iii) of \ref{C1} we conclude that $\Phi(\{t\} \times J)$ is an isogeny leaf.

Conversely, let $I' \subset I(y) \subset W'$ be an isogeny leaf.
Choose an irreducible component $J' \subset \rho^{-1}(I')$ and $w' \in J'(k)$ with  $ \rho(w') = y$. We have $\pi(J') =: z' \in C$, and we can choose  $(t_0,u_0) \in T \times J$ with $\Psi(t_0,u_0) = w'$, using that fact that $\Psi$ is surjective. As $\pi(\Psi(\{t\} \times J)) = z'$, and because $\Psi$ is proper and finite, and $\dim(J') = \dim(J)$ we conclude that 
$\Psi(\{t\} \times J) = J'$. Hence $\Phi(\{t\} \times J) = I'$. This shows (iv). This finishes the proof of  \ref{C2}, and the proof of Theorem \ref{TI} is concluded.     \B\ref{TI}

\subsection{}\label{cap}{\bf Corollary.}  {\it Let $W''$ be an irreducible component of 
$\cW^0_{\xi}(\cA_g \otimes k)$.}\\
{\bf (i)} {\it Let $I \subset W''$ be an isogeny leaf and $C \subset W''$ a central leaf.
Then $C \cap I \not= \emptyset$.}\\
{\bf (ii)} {\it Every isogeny leaf in $W''$ has dimension
$\dim(W'') - c(\xi)$.
}\\
Indeed, if $I = \Phi(\{t\} \times J)$ and $C = \Phi(T \times \{u\})$, then 
$\Phi(t,u) \in C \cap I$. 

As $\Phi$ is finite and surjective we obtain $\dim(W'') = \dim(T) + \dim(J)$;  for any isogeny leaf $I$ we have $\dim(J) = \dim(\Phi(\{t\} \times J)) = \dim(I)$. \B

\subsection{}\label{ixi}{\bf Corollary / Notation}  $i(-)$. {\it Let $\xi$ 
be a symmetric Newton polygon 
(dimension $g$, height $2g$). There exists an integer $i(\xi)$ such that for every
$x \in \cA = \cA_{g,1} \otimes \FF_p$ the dimension of every geometric irreducible 
component of
$I(x) \subset \cA$ equals $i(\xi)$. We have $c(\xi) + i(\xi)  = \sdim(\xi)$.}
\\
Reminder: The integer $\sdim(\xi)$ was defined in \cite{CH}, 3.3, and 
in \cite{NP}, 1.9; in \cite{NP}, Theorem 4.1  
we have proved that
the dimension of every geometric component of $\cW_{\xi}(\cA_{g,1} \otimes \FF_p)$ 
equals this integer: 
$$\dim\left(\cW_{\xi}(\cA_{g,1} \otimes \FF_p)\right) \quad=\quad \sdim(\xi).$$  \B

\subsection{}
Note that in the last corollary we consider {\it principally polarized} abelian varieties. 
{\it The dimension of a geometric component of an isogeny leaf depends 
on the component of the Newton polygon stratum containing the leave}. \\
{\bf Example.} For $\sigma_3 = \sigma = \xi = 3{\cdot}(1,1)$ we have $i(\xi)=2$, 
i.e. isogeny leaves in the 
supersingular, principally polarized $g=3$  locus all have dimension equal to two;
note that any supersingular central leaf is finite;
however there do exist supersingular isogeny leaves in 
$\cA_{g,d} \otimes \FF_p$ of dimension three,  see \cite{K.FO}, 6.10 ($d = p^3$), 
\cite{Li.FO}, 12.4  ($d = p^2$);
this produces isogeny leaves in $\cW_{\sigma}(\cA_3 \otimes \FF_p)$ of dimension 
bigger than $2 = i(\sigma_3)$.

Also for some non-supersingular Newton polygons  analogous examples can be given, showing that dimensions of isogeny leaves related to one Newton polygon, contained in different components of $\cW^0_{\xi}(\cA_g \otimes k)$, need not be the same. Compare \ref{maxdim}.

\subsection{}{\bf Corollary.} {\it Let $I \subset \cW_{\xi}(\cA_g \otimes k)$ be an isogeny leaf. Then:\\
(i) there is a unique irreducible component $W''$ of $\cW_{\xi}(\cA_g \otimes k)$ containing $I$;\\
(ii) there is a dense open subset $U \subset I$ such that for every $y \in U$ we have $I(y) = I$.
}
\\
{\bf Proof.} Suppose $I \subset W''$, an irreducible component  $\cW_{\xi}(\cA_g \otimes k)$.    We consider $w \in W(k)$ as in  the proof of Theorem \ref{TI} in particular  such that its image in $W'' \subset \cW_{\xi}(\cA_g \otimes k)$ is contained in no other irreducible component of $\cW_{\xi}(\cA_g \otimes k)$. Write $x'' = \rho''(w)$, and let $C_{x''} \subset W''$ the central leaf passing through this point. By the choice of $w$ and by (i) of Theorem \ref{dim=} it follows that $C_{x''}$ does not meet any component of  $\cW_{\xi}(\cA_g \otimes k)$ different from $W''$. By \ref{cap} we see that $I \cap C_{x''}$ is non-empty; hence $I$ contains at least one point contained in no other component but $W''$.  This proves (i). 

We consider in the proof of Theorem \ref{TI} in particular and moreover a choice of $w \in W(k)$ such that $W'$ is normal at $x' = \rho(w)$, i.e. such that $W' \leftarrow W$ is birational at $w$. Note that $(W')^{/x'} \cong (W')^{/y}$ for every $y \in C_{x'}(k)$;  hence  $W'$ is normal at all points of $C_{x'}$. This means that $W' \leftarrow W$  is birational at all points above $C_{x'}$.

Let $I'_1, \cdots ,I'_s$ be the irreducible components of the inverse image of $I$ under $W'' \leftarrow W'$. 
By the previous argument, and by \ref{cap} we conclude that every such component meets the normal locus of $W'$. Hence the inverse image of $I$ under $W'' \leftarrow W$ equals $I_1 \cup \cdots \cup I_s$, union of the inverse images. For every index $i$ consider $U_i \subset I_i$ as the set where $I_i$ does not meet other irreducible components of the corresponding fiber of $W \to C$. Let $V'' \subset I$ be the set of points, images of non-normal points of $W'$  intersected with $I$; construct $U$ as the intersection of $I - V''$ with all images $\rho''(U_i)$. For $y \in U(k)$ an irreducible component of $I(y)$, i.e. an isogeny leaf passing through $y$, is an image of some $\{t\} \times J$, as we know by Theorem \ref{TI}, hence comes from an irreducible component of the fiber $W \to C$ above a point $\pi(I'_i)$; by the choice of $U$ this means that $I(y)$ is irreducible at $y$, which proves the corollary.  \B

\subsection{}{\bf Remark.} Although an isogeny leaf is contained in a unique component of its Newton polygon stratum, a central leaf may be contained in more that one component of its Newton polygon stratum.

\subsection{}{\bf Remark.} Using \cite{C.FO}, see \ref{cxi}, we  conclude that for 
symmetric Newton polygons
$\xi \succneqq \xi'$ we have $c(\xi) > c(\xi')$ and $i(\xi) \leq i(\xi')$. 

\subsection{}{\bf Remark.} Isogeny leaves,  in the context of 
Rapoport-Zink spaces as in \cite{RZ}, and the product structure 
as given in \ref{TI} are studied and described in the case of certain Shimura varieties 
in \cite{EM}. It would be worthwhile to describe the product structure in general for 
all Shimura varieties in mixed characteristic, in particular for the Siegel 
case studied here.

\section{Some questions}
\subsection{}\label{HO}{\bf Conjecture} (the Hecke Orbit Conjecture). Consider a point 
$[(A,\lambda)] = x \in \cA = \cA_g \otimes \FF_p$, 
and consider its 
Hecke orbit $\cH(x) \subset \cA$. {\it We expect this Hecke orbit to be dense in 
its Newton polygon stratum in the moduli space, i.e. the Zariski closure is expected to be:}
$$\overline{\cH(x)} \quad\stackrel{?}{=}\quad \cW_{\cN(A)}(\cA).$$
Notation: $\cH(x)$ is the set of all $y = [(B,\mu)]$ such that there exists a field $L$,
an isogeny $\va: A_L \to B_L$, and an integer $n \in \ZZ_{>0}$ such that $\va^{\ast}(\mu_L) = n{\cdot}\lambda_L$.  We write $y \in \cH_{\ell}(x)$
if moreover $\ell$ is a prime number and the degree of $\va$ and $n$   are powers of $\ell$. 

\vn
Note that  Hecke-$\ell$-actions ``move'' in a central leaf: under $\ell$-degree-isogenies, 
with $\ell \not= p$ the $p$-divisible groups are unchanged. Note that  Hecke-$\alpha$-actions 
(isogenies with local-local kernel) ``move''  in an isogeny leaf. Moreover we have proved that a non-supersingular Newton polygon in $\cA_{g,1} \otimes \FF_p$ is geometrically irreducible, see \cite{NPirred}.
Therefore  Conjecture \ref{HO} would follow, using \ref{TI}, in case we can prove: 

\subsection{}\label{HaO}{\bf Conjecture.} {\it For every prime number $\ell$, different
from $p$, the Hecke-$\ell$-orbit   $\cH_{\ell}(x)$ in $\cA = \cA_{g,1} \otimes \FF_p$ 
is Zariski-dense in 
$\cC(x)$, a finite union of central leaves.}

These conjectured results seem to be true; proofs still have to be written out in full.

\subsection{} Let $\cX \to S$ be a $p$-divisible group over a   scheme over $k$.
Define a reduced subscheme $I \subset S$ to be $H_{\alpha}$ if there exists a chart
$(T \twoheadrightarrow I, \va: Y \times T \to \cX)$, as in \ref{defH}, where $T \to I$ 
is a morphism of finite type. Let $x \in S(k)$.\\
{\bf Question.} {\it Does there exist a maximal $I(x) \subset S$ subscheme, as in \ref{isogl},
which is} $H_{\alpha}$ ?

\subsection{}\label{maxdim}{\bf Conjecture.} {\it Let $\xi$  be a symmetric Newton polygon, with $p$-rank 
equal to $f$, 
i.e. $\xi$ has exactly $f$ slopes equal to zero.  We expect that $\cW_{\xi}(\cA_g)$ 
has a component 
of dimension precisely $(g(g-1)/2) + f$.} 

Note that it is clear that every such component
has at  most this dimension, see \cite{N.O}.

\subsection{}\label{ConjCx}{\bf Conjecture.} Let $\xi$ be a symmetric Newton polygon, of height $2g$,
not equal to the supersingular one: $\xi \succneqq \sigma$. {\it We expect that 
$\cW_{\xi}(\cA_{g,1} \otimes \FF_p)$ is}  geometrically irreducible. {\it We expect that for every non-supersingular $x \in  \cA_{g,1}(k)$  we have $C(x) = C_x$; i.e. we expect that $C(x)$ is} (geometrically) irreducible.

These conjectured results seem to be true; proofs still have to be written out in full.

\subsection{}{\bf Expectation.} Let $\xi$ be a symmetric Newton polygon, of height $2g$,
and let $\eta$ be the generic point of a component of $\cW_{\xi}(\cA_g \otimes \FF_p)$;
let $X_{\eta}$ be the $p$-divisible group over $k = \overline{\FF_p(\eta)}$ derived
from the geometric fiber over $\eta$. By \cite{Z2}, especially Lemma 9 and Proposition 12, 
there exists a (canonical, minimal) choice for an isogeny $Z \to X_{\eta}$ such that
$Z$ is \csd. {\it We expect that} $Z_k \cong_k H(\xi)$; here $H(\xi)$ is the minimal $p$-divisible group with Newton polygon $\xi$, see \ref{min}. 

\subsection{}\label{i} Let us choose a number $i \in \ZZ_{>0}$. For any point 
$[(X,\lambda)] = x \in \cA_g \otimes \FF_p$ we can consider $\va:= (X,\lambda)[p^i]$,
and we can study
$S^{(i)}_{\va}(\cA_g \otimes \FF_p)$, the set of points $y = [(Y,\mu)]$ such that there
exist: an algebraically closed field $k$ over which $y$ is defined, and  an isomorphism
$(X,\lambda)[p^i] \otimes k \cong (Y,\mu)[p^i] \otimes k$; probably this is a locally 
closed set in $\cA_g \otimes k$.

\vn
Choosing $i=1$ we obtain $S^{(1)}_{\va}(\cA_g \otimes \FF_p) = S_{\va}$, 
the EO-strata as defined in \cite{EO}.\\
Note that the leaves defined by $S^{(i+1)}$ are contained in leaves defined by $S^{(i)}$: for
$\va_1 = , \va_2, \cdots, \va_i = [(X,\lambda)[p^i]], \cdots $ all coming from the same $(X,\lambda)$ 
we obtain 
$S^{(1)}_{\va_1}(-) \supset  S^{(2)}_{\va_2}(-) \supset \cdots$; this descending chain stabilizes 
after a finite number of steps.\\
For $i>>1$ we obtain central leaves: given $g$, there exists $N$ such that for every $x$ we have
$S^{(N)}_{\va_N}(-) = \cC_x(-)$, see \ref{2}.

\vn
We studied $S^{(1)}$ in \cite{EO}, and we consider $S^{(N)} = S^{({\infty})}$ in this paper;
one could also study the ``intermediate'' cases $S^{(i)}$. 

\subsection{}{\bf Conjecture.} Let $(H(\xi),\zeta)$ be a principally quasi-polarized
minimal $p$-divisible group, see \ref{min}. We expect: {\it the central stream, i.e. the central leaf 
$\cZ_{\xi} = \cC_{H(\xi)}(\cA_{g,1} \otimes \FF_p)$, is the} EO-{\it stratum} $S^{(1)}_{H(\xi)[p]}$ 
{\it  associated with $H(\xi)[p]$;
conversely, every central leaf equal to its} EO-{\it stratum is a central stream.} See \cite{Min}.

\subsection{}{\bf Conjecture.} Let $\xi$ be a symmetric Newton polygon, 
let $\va = \va(\xi)$ be the isomorphism class of $H(\xi)[p]$. Let $\psi$ be the  isomorphism 
class
of a BT$_1$ such that $S_{\psi}(\cW^0_{\xi}(\cA_{g,1} \otimes \FF_p)) \not= \emptyset$.
We expect that in case: $\va(\xi) =\va \subset \psi$, in the terminology of \cite{EO}, (14.3); this means: $S_{\va}$ is contained in the Zariski closure  of $S_{\psi}$.

\subsection{}{\bf Question.} We try to study the closure of central leaves in lower 
Newton polygon strata. Let $\xi$ be a symmetric Newton polygon, $\cZ_{\xi}$ the 
central stream, and $C \subset W^0_{\xi}$ be a central leaf. It might be  that $\cZ_{\xi}$
and $C$ have the same boundary inside $\cA_g$:
$$\left(\cZ_{\xi}^c - \cZ_{\xi}\right) \cap \cA_{g,1} \quad\stackrel{?}{=}\quad
\left(C^c - C\right) \cap \cA_{g,1}.$$
Moreover it could be  that for $\xi' \prec \xi$ we have $\cZ_{\xi'} \subset (\cZ_{\xi})^c$; 
note that we expect in
general that {\it a leaf in $W^0_{\xi'}$ is not in the closure of a leaf in $W^0_{\xi}$}
(depending on $\xi' \prec \xi$ and on which leaf is chosen).

\subsection{} We can also study Hecke-correspondences which have only local-\'etale 
and \'etale-local kernels. Orbits under such correspondences we denote by $\cH_{pna}(-)$, where
$p$ stands for of $p$-power degree, and $na$ for non-$\alpha_p$. Clearly such correspondences 
do not change the $p$-divisible group up to isomorphism over $k$. I have no reasonable guess 
whether an orbit under $\cH_{pna}$ could give a dense subset in a central leaf. 
  
%\newpage

\section{Notations}\label{Not}
\n
\subsection{} In this paper all base rings, and all base schemes will
 be in characteristic $p$. Usually we will write $K$ for an arbitrary 
 field (in characteristic $p$), 
and $k$ for an algebraically closed field.

\subsection{}\label{N} We say that a scheme $S$ satisfies condition (N) if $S$ is integral, and the normalization $S' \to S$ gives a noetherian scheme  $S'$.

Examples of such schemes: an integral  scheme of finite type over a field; an integral Japanese scheme, see \cite{EGA}, IV$^1$, Chap. 0, 23.1.1; an integral excellent scheme, see \cite{EGA}, IV$^2$.7.8.

\subsection{}
We write $\cA_g \to \ZZ$ for the moduli scheme of polarized abelian schemes of 
relative dimension $g$. 
We write $\cA$ for the moduli spaces of  principally polarized 
abelian varieties in characteristic $p$:
$$\cA = \cA_{g,1} \otimes \FF_p;$$ 
in some considerations
all degrees of polarizations will be allowed.
In case
we also consider a level structure, we will only consider
a level-$n$-structure with $n \in  \ZZ_{>0}$ and prime to $p$. 
We write $\cA_{g,\ast} = \cup_d \ \ \cA_{g,d}$  and $\cA_{g,\ast,\ast} = 
\cup_{d,n}\ \  \cA_{g,d,n}$.

\subsection{}\label{NP}{\bf Newton polygons.} {\it The set of isogeny classes of $p$-divisible groups 
over any algebraically closed field is the same as the set of Newton polygons.} 
This combinatorial, discrete invariant can be described as follows.
For given integers $d \in \ZZ_{\geq 0}$ 
and $h \in \ZZ_{\geq d}$, writing 
$h = d+c$,  we will understand under a {\it Newton polygon belonging to 
$d$ and $h$:} 

a polygon starting at $(0,0)$, ending at  $(h,c)$, 

which is lower convex, 

has breakpoints 
in $\ZZ \times \ZZ$, 

and which has slopes $0 \leq \lambda \leq 1$.\\
Note that the set of all $p$-divisible groups over a given 
algebraically closed field $k$ up to isogeny is in 
bijective correspondence with the set of all Newton polygons, 
as was shown in the Dieudonn\'e - Manin theory, see \cite{Manin}. 
Let us make this correspondence precise.  For given non-negative 
coprime integers $m$ and $n$, we fix a $p$-divisible group $G_{m,n}$
as in \cite{Manin}; this is defined over $\FF_p$; 
this $p$-divisible group has dimension $m$ and its Serre-dual 
$(G_{m,n})^t = G_{n,m}$ has dimension $n$. In fact, $G_{1,0} = \GG_m[p^{\infty}]$,
and $G_{0,1} = \underline{\QQ_p/\ZZ_p}$, and for positive $m$ and $n$ the covariant
Dieudonn\'e module $\DD(G_{m,n})$ is generated by one element $e$ satisfying the relation
$\cF^m{\cdot}e = \cV^n{\cdot}e$. To the $p$-divisible group $G_{m,n}$ we associate the
Newton polygon $\cN(G_{m,n})$ consisting of $m+n$ slopes equal to $n/(m+n)$;
$$\dim(G_{m,n}) = m, \ \  \dim((G_{m,n})^t) = n, \ \ \  
\mbox{\rm slope}(G_{m,n}) = \frac{n}{m+n}, \ \ \mbox{\rm height}(G_{m,n}) = h = m+n.  $$
For any $p$-divisible group $X$ over a field $K$, there exists an isogeny
$$X \otimes k \quad \sim \quad \bigoplus_i \ \ G_{m_i,n_i}$$
over an algebraically closed field $k$ containing $K$; 
the slopes $n_i/(m_i+n_i)$ with multiplicity $m_i+n_i$ appearing in this sum are  arranged in 
non-decreasing order and they 
give the Newton polygon of $X$. \\
{\bf Remark.} Slopes are defined as e.g. in \cite{Katz}, \cite{CH}, \cite{NP}. This means: 
for any slope $\lambda'$ appearing \cite{Z2}, and in \cite{FO.Z}, the 
slope defined here is equal to $\lambda = 1 - \lambda'$. In definitions adapted from those
papers in this paper we will carry out this adaptation. In this paper a slope filtration,
the natural slope filtration, and a Newton polygon, is according to non-decreasing slopes:
$1 \leq \lambda_1 \leq \cdots \leq \lambda_i \leq \cdots \leq 1$.\\
{\bf Remark.} For a scheme $T$ over a scheme $S$  we will write $F_{T/S}= F: T \to T^{(p)}$
for the relative Frobenius, and for a commutative group scheme we write 
$V_{G/S} = V: G^{(p)} \to G$ for the relative Verschiebung. Whenever covariant
Dieudonn\'e module theory is considered we write $\cV$ and $\cF$ for the Verschiebung and
Frobenius maps; note that in the covariant theory $\DD(F) = \cV$, 
and $\DD(V) = \cF$ (in the obvious sense); see \cite{CH}, Section 1 for an explanation.

\vn
Note that we use ``$V$-slopes'' on $p$-divisible groups and abelian varieties, which is the same as ``$\cF$-slopes'' on Dieudonn\'e modules.

\vn
A Newton polygon is called {\it symmetric} if any slope  $\lambda$ is appearing with the same 
multiplicity as $1-\lambda$. The Newton polygon of a $p$-divisible group $X$ is
symmetric if and only if $X$ and its Serre-dual $X^t$ are isogenous. The Newton polygon of
(the $p$-divisible group of) an abelian variety $A$ is symmetric: a polarization on an 
abelian variety by the duality theorem, see \cite{CGS}, 18.1, gives an 
isogeny between $A[p^{\infty}] =: X$ and $X^t$. Note that $A^t[p^{\infty}] = (A[p^{\infty}])^t$ by the duality theorem, see \cite{CGS}, 19.1.

\vn
We write $\gamma \prec \beta$ if these are Newton polygons related with the same numbers
$d$ and $c$ and no point of $\gamma$ is below $\beta$ (and we will say  ``$\gamma$
is lying above $\beta$'').

\vn
Let $\cX \to S$ be a family of $p$-divisible groups, and let $\beta$ be a Newton polygon. We write $\cW^0(S)$ for the set of points $s \in S$ such that the Newton polygon of $\cX_s$ equals $\beta$. We write $\cW(S)$ for the set of points $s \in S$ such that the Newton polygon of $\cX_s$  is lying above $\beta$. By Grothendieck-Katz, see  \cite{Katz}, Th. 2.3.1, we know that  $\cW(S) \subset S$ is closed, and $\cW^0(S) \subset S$ is locally closed.

\subsection{}\label{min}
{\bf Minimal $p$-divisible groups.} For a pair $(m,n)$ of coprime non-negative integers 
we write $H_{m,n}$ for the $p$-divisible group as in \cite{Purity}, 5.3; $H_{m,n}$ can be 
defined over $\FF_p$, it is  isogenous with $G_{m,n}$ and over
$\overline{\FF_p}$ its endomorphism ring is the maximal order in its endomorphism algebra;
these properties characterize $H_{m,n}$ over an algebraically closed field;
also see \cite{Min}.
Note that $\End(H_{m,n} \otimes \FF_p)$ is a commutative algebra,
free of rank $h=m+n$ over $\ZZ_p = W(\FF_p)$; in this ring we have a uniformizer 
$\pi$ and relations  $\pi^h=p$, and  $\pi^m = \cF$, and $\pi^n = \cV$;  for any field 
$K \supset \FF_{p^h}$ the
algebra  $\End(H_{m,n} \otimes K)$ is not commutative if $h=m+n>1$,  it is free 
of rank $h$
over $W(\FF_{p^h})$, it is maximal in its full ring of quotients, and it  has $\pi$ 
as uniformizer.

For a Newton polygon $\beta = \oplus (m_i,n_i)$ we write
$H(\beta) := \oplus H_{m_i,n_i}$; this will be called the minimal $p$-divisible group
with Newton polygon equal to $\beta$. For any field $K \supset \FF_p$ we write
$H(\beta)$ instead of $H(\beta) \otimes_{\FF_p} K$, if no confusion can occur. 

\subsection{}\label{defcsd}{\bf Completely slope divisible.}
We take the definition, due to T. Zink,  of a completely slope divisible $p$-divisible 
group from
\cite{Z2}, \cite{FO.Z}. This means:\\
{\bf Definition.} {\it Let $X \to S$ be a $p$-divisible group.  
Let $0 \leq t_1 < t_2 < \ldots < 
t_m \leq s$ be integers and equalities. A $p$-divisible group 
$Y$ over a scheme $S$ is said to be} {\rm completely slope 
divisible with respect to} {\it these integers if $Y$ has a filtration
 by $p$-divisible subgroups:
\begin{displaymath}
   0 = Y_0 \subset Y_1 \subset \ldots \subset Y_m = Y
\end{displaymath}
such that the following properties hold:}\\
(1) {\it The quasi-isogenies
\begin{displaymath}
      \frac{F^s}{p^{s-t_i}} : Y_i \rightarrow Y_i^{(p^s)}
\end{displaymath} 
are isogenies for $i = 1, \ldots ,m$.}\\
(2) {\it The induced morphisms:
\begin{displaymath}
   \frac{F^s}{p^{s-t_i}} : Y_i/Y_{i-1} \rightarrow (Y_i/Y_{i-1})^{(p^s)}
\end{displaymath}
are isomorphisms.}

\n
Note that the last condition implies that $Y_i/Y_{i-1}$ is
isoclinic of exact slope $\lambda_i := t_i/s$. See \cite{FO.Z}, 1.5 for
a characterization of a $p$-divisible group over a field being \csd.
Note that we do not suppose that $t_i$ and $s$ are relatively prime; 
it is easy to give an example of a $p$-divisible group over a field
which is \csd \  with respect to $d{\cdot}t_i$ and $d{\cdot}s$,
where $d\in \ZZ_{>1}$ and not \csd \  with respect to $t_i$ and $s$.

\subsection{} We write $\va: (A,\lambda) \to (B,\mu)$ in case $(A,\lambda)$ and
$(B,\mu)$ are polarized abelian schemes over some base $S$, and $\va: A \to B$
is an isogeny such that $\va^{\ast}(\mu) = (\lambda)$. We write 
$\va: (A,\lambda,f) \to (B,\mu,f)$ in case moreover $f$ is a level structure
prime to the characteristic on the base, and prime to the degree of $\va$. We
use analogous notations for $p$-divisible groups.

\section{Some terminology and an example}\label{T}
We use: $K$ for a field of characteristic $p$, all base schemes will be in characteristic $p$,

we write $k$ for an algebraically closed field 
of characteristic $p$;

we write $A$ for an abelian variety, $X$ for a $p$-divisible group over a field, 

and usually 
$\cX$ for a $p$-divisible group over a base of positive dimension;

we write $G$ for a finite group scheme, 

and usually 
$\cG$ for a  group scheme finite and locally free  over a base of positive dimension.

\vn
An isoclinic Newton polygon has all slopes equal to each other; the only symmetric 
isoclinic Newton polygon 
is the supersingular one: $\sigma = g{\cdot}(1,1)$. 

\vn
In this paper we have defined the notion of {it central leaf}, see \ref{cl},  
\ref{pcl} and \ref{defcl}. W used this terminology ``{\it central leaf}'' 
in order to contrast with ``{\it isogeny leaves }''. If no isogeny leaves 
are in consideration in a future publication one can simply 
say ``leaf'' and ``foliation'' for the notion of central leaves used here.

We think that the central leaf (leaves) connected with a minimal $p$-divisible group, 
see \ref{min}, will play an important role; hence we gave this a name, the ``central stream'',
in order to distinguish this from the other central leaves.

\vn
For ``central leaf'', $C(x)$ and $C_{x}$ see:  \ref{defcl}; \ \ for ``isogeny leaf'', and  $I(x)$ : \ref{defil};\\
``central stream'': \ref{defcs};  
\ \ ``\gfc'': \ref{defgfc}; \ \ ``constant'': \ref{defconst};\\ 
``completely slope divisible'':  \ref{defcsd}; 
\ \ ``minimal'': \ref{min};\\
$\dim(\beta)$:  \cite{NP}, 1.6; \quad  $\sdim(\xi)$: \cite{NP}, 1.9;\\ 
$\cu(\beta)$: \ref{cu}; \quad $\c(\xi)$: \ref{dim=};  \quad $i(\xi)$ : \ref{ixi}.

\subsection{}\label{exa} As illustration we record for $g=4$ the various data considered:
$$
\begin{array}{|c|l|c|c|c|c|c|}\hline
{\rm NP} &\xi& f &sdim(\xi)& c(\xi)& i(\xi)& {\rm ES}(H(\xi))\\
\hline
\rho&(4,0)+(0,4)&4&10&10&0&(1,2,3,4)\\  \hline
f=3&(3,0)+(1,1)+(0,3)&3&9&9&0&(1,2,3,3)\\  \hline
f=2&(2,0)+(2,2)+(0,2)&2&8&7&1&(1,2,2,2)\\ \hline
\beta&(1,0)+(2,1)+(1,2)+(0,1)&1&7&6&1&(1,1,2,2)\\  \hline
\gamma&(1,0)+(3,3)+(0,1)&1&6&4&2&(1,1,1,1)\\  \hline
\delta&(3,1)+(1,3)&0&6&5&1&(0,1,2,2)\\  \hline
\nu&(2,1)+(1,1)+(1,2)&0&5&3&2&(0,1,1,1)\\  \hline
\sigma&(4,4)&0&4&0&4&(0,0,0,0)\\  \hline

\end{array}
$$
Here $\rho \succ (f=3) \succ (f=2) \succ \beta \succ 
\gamma \succ \nu \succ \sigma $ and $\beta \succ \delta \succ\nu$. 
The notation ES, encoding the isomorphism type of a BT$_1$ group scheme, 
is as in \cite{EO}; the number $f$ indicates the $p$-rank.

\noindent
\begin{tabbing}
\hspace{80 mm}
\=Frans Oort \\ 
\>Mathematisch Instituut         \\
 \>Postbus 80.010                     \\   
  \>NL-3508 TA Utrecht  \\ 
\>The Netherlands      \\ 
\>  \ \ \ \ \ \ email:   oort@math.uu.nl       \\
\end{tabbing}

%\RRR

\newpage
\begin{thebibliography}{}
%\bibitem{Chai}
%C.-L. Chai -- {\it  Every ordinary symplectic isogeny class in 
%positive characteristic is dense in the moduli.}
%Invent. Math. {\bf 121} (1995), 439-479. 

\bibitem{C.FO}
C.-L. Chai \& F. Oort -- {\it  Canonical coordinates on leaves of $p$-divisible groups.} [In preparation]

\bibitem{FC}
G. Faltings \& C.-L. Chai --  {\it   Degeneration of abelian varieties.}
Ergebnisse Bd 22, Springer - Verlag, 1990.

\bibitem{EGA}
A. Grothendieck \& J. Dieudonn\'e --  {\it Elements de g\'eom\'etrie alg\'ebriques.}\\ 
%Ch. III$^1$: {\it \'Etude cohomologique des faisceaux %coh\'erents.} 
%Publ. Math. IHES {\bf 11} (1961).\\
Ch. 0 (suite): IV$^1$, Publ. Math. IHES {\bf 20} (1964).\\
Ch. IV$^2$: {\it Etude locale des sch\'emas et des morphismes de sch\'emas.}
Publ. Math. IHES {\bf 24} (1965).

\bibitem{HT}
M. Harris \& R. Taylor --  {\it  The geometry and cohomology of some simple Shimura 
varieties.} Ann. of Math. Studies 151, Princeton Univ. Press 2001.

\bibitem{HAG}
R. Hartshorne  --  {\it Algebraic Geometry.} Grad. Texts 52; Springer - Verlag 1977.


\bibitem{dJ}
A. J. de Jong  --  {\it Crystalline Dieudonn\'e module theory via
formal rigid geometry.}  Publ. Math. IHES {\bf 82} (1995), 5-96.

\bibitem{Purity} 
A. J. de Jong \& F. Oort  --   
{\it  Purity of the stratification by Newton polygons.}  
Journ. A.M.S. {\bf 13} (2000), 209-241.

\bibitem{K.FO}
T. Katsura \& F. Oort  --  
{\it  Supersingular abelian varieties of dimension 
two or three and   class numbers.}  Algebraic Geometry, 
Sendai 1985 (Ed. T. Oda). Adv. Stud. in Pure Math.      
{\bf 10} (1987), Kinokuniya Cy Tokyo and North-Holland 
Cy Amsterdam, 1987  ; pp.  253-281. 

\bibitem{KatzST}
N. M. Katz  --  {\it  Serre-Tate local moduli.}
In: Surfaces alg\'ebriques (S\'em. de g\'eom. alg\'ebr.
d'Orsay 1976-78). Lecture Notes Math. 868, Springer - Verlag 1981; Exp. V-bis, 
pp. 138-202.

\bibitem{Katz}
N. M. Katz -- {\it Slope filtration of
$F$-crystals.} Journ. G\'eom. Alg. Rennes, Vol. I, 
Ast\'erisque  {\bf 63} (1979),
Soc. Math. France,    113-164. 




\bibitem{Li.FO}
K.-Z. Li \& F. Oort -- {\it Moduli 
of supersingular abelian varieties.}  Lecture Notes Math. 1680, 
Springer  1998; 116 pp.



\bibitem{Manin}
Yu. I. Manin -- {\it The theory of 
commutative formal groups 
over fields of
finite characteristic.} Usp. Math. 
{\bf 18} (1963), 3-90; 
Russ. Math. Surveys {\bf 18}
(1963), 1-80.

\bibitem{EM}
E. Mantovan --  {\it  On certain unitary group Shimura varieties.} 
Harvard PhD-thesis, April 2002. [To appear]

\bibitem{N.O}
P. Norman \& F. Oort  -- 
{\it Moduli of abelian varieties.} Ann.  Math.   
{\bf 112} (1980), 413-439.

\bibitem{Oda}
T. Oda --  {\it The first de Rham cohomology group and Dieudonn\'e modules.} 
Ann. Scient. Ec. Norm. Sup. {\bf 2} (1969), 63-135.

\bibitem{CGS}
F. Oort --  {\it Commutative group schemes.} 
Lecture Notes Math. 15, Springer - Verlag 1966.  

\bibitem{CH}
F. Oort --  {\it  Newton polygons and formal groups: conjectures by Manin and Grothendieck.}  Ann. Math. {\bf 152} (2000), 183-206.


\bibitem{EO}
F. Oort  --  {\it A stratification of a moduli space of 
polarized abelian varieties.} In: {\it Moduli of abelian varieties.}  
(Ed. C. Faber, G. van der Geer, F. Oort). Progress Math. 195, 
Birkh\"auser Verlag 2001; pp. 345-416.



\bibitem{NP}
F. Oort -- {\it Newton polygon strata in the 
moduli space of abelian varieties.} In:   
{\it Moduli of abelian varieties.}  (Ed. C. Faber, 
G. van der Geer, F. Oort). Progress Math. 195, Birkh\"auser 
Verlag 2001; pp. 417 - 440. 

\bibitem{Min}
F. Oort  --  {\it Minimal $p$-divisible groups.} [To appear] 

\bibitem{NPirred}
F. Oort -- {\it  Irreducibility of Newton polygon strata.} [In preparation]


\bibitem{FO.Z}
F. Oort \& Th. Zink --  {\it   Families of $p$-divisible groups with 
constant Newton polygon.} Documenta Mathematica {\bf 7} (2002), 183-201, see\\
http://www.mathematik.uni-bielefeld.de/ documenta/vol-07/09.html


\bibitem{RZ}
M. Rapoport \& Th. Zink --  {\it   Period spaces for $p$-divisible groups.}
Ann. of Math. Studies 141, Princeton Univ. Press 1996.

\bibitem{Reiner}
I. Reiner --  {\it Maximal orders.} Academic Press, 1975.

\bibitem{Z2}
Th. Zink -- {\it On the slope filtration.} 
Duke Math. J. Vol. {\bf 109} (2001),   79-95.
	
	
	

 
 \end{thebibliography}
\end{document}